\documentclass[aos,seceqn,nameyear,dvips]{arximspdf}
\usepackage{graphics}
\doi{10.1214/09-AOS765}
\volume{38}
\issue{3}
\pubyear{2010}
\firstpage{1733}
\lastpage{1766}

\makeatletter

\newproclaim{assumption}{Assumption}[section]
\newtheorem{proposition}{Proposition}[section]

\newproclaim{example}{Example}[section]
\newtheorem{corollary}{Corollary}[section]
\newtheorem{lemma}{Lemma}[section]

\makeatother

\begin{document}
\begin{frontmatter}

\title{Approximation of conditional densities by smooth mixtures of
regressions}
\runtitle{Approximation by mixtures of regressions}

\begin{aug}
\author[A]{\fnms{Andriy}
\snm{Norets}\corref{}\ead[label=e1]{anorets@princeton.edu}\ead[label=u1,url]{http://www.princeton.edu/\texttildelow anorets}}
\runauthor{A. Norets}
\affiliation{Princeton University}
\address[A]{313 Fisher Hall\\
Department of Economics\\
Princeton University\\
Princeton, New Jersey 08544\\
USA\\
\printead{e1}} %adresu isvedimo komanda gale!
\end{aug}

% HISTORY:
\received{\smonth{8} \syear{2009}}
\revised{\smonth{11} \syear{2009}}

% ABSTRACT
%
\begin{abstract}
This paper shows that large nonparametric classes of conditional
multivariate densities can be approximated in the Kullback--Leibler
distance by different specifications of finite mixtures of normal
regressions in which normal means and variances and mixing
probabilities can depend on variables in the conditioning set
(covariates). These models are a special case of models known as
``mixtures of experts'' in statistics and computer science literature.
Flexible specifications include models in which only mixing
probabilities, modeled by multinomial logit, depend on the covariates
and, in the univariate case, models in which only means of the mixed
normals depend flexibly on the covariates. Modeling the variance of the
mixed normals by flexible functions of the covariates can weaken
restrictions on the class of the approximable densities. Obtained
results can be generalized to mixtures of general location scale
densities. Rates of convergence and easy to interpret bounds are also
obtained for different model specifications. These approximation
results can be useful for proving consistency of Bayesian and maximum
likelihood density estimators based on these models. The results also
have interesting implications for applied researchers.

%Obtained characterizations of approximable classes of densities for
%different specifications of the model can be useful for choosing
%appropriate specifications in applications.
\end{abstract}

% KEYWORDS
%

\end{frontmatter}

%s1 ###
\section{Introduction}

This paper explores approximation properties of finite smooth mixtures
of normal regressions as flexible models for conditional densities.
These models are a special case of
mixtures of experts (ME) introduced by \citet{JacobsEtAl91}.
ME have become increasingly popular is statistical literature since
they are very flexible, easy to interpret and
reasonably easy to estimate. See, for example, papers by \citet
{JordanJacobs94}
and \citet{JordanXu95} who employ
the expectation maximization (EM) estimation algorithm
or papers by
\citet{PengJacobsTanner1996},
\citet{WoodJiangTanner02}, \citet{Geweke07} and\break \citet
{VillaniKohnGiordani07} who use Markov chain Monte Carlo methods for
estimation of ME in the Bayesian framework.
This paper contributes to the literature that provides a theoretical
explanation of the success of ME models in applications.
In particular, I show that large classes of conditional densities can
be approximated in the Kullback--Leibler (KL) distance by finite smooth
mixtures of normal regressions.
Approximation results are obtained in the KL distance for the following
reason. If a data generating density is in the KL closure of a class of
models then this density can be consistently estimated from data by
these models under weak regularity conditions [see, e.g.,
\citet{GhoshRamamoorthi03} for a textbook treatment of Schwarz's
theorem on posterior consistency and
\citet{RoederWasserman97} for posterior consistency results for
finite mixture of normals].

Consider a joint probability distribution $F$ on a product space $Y
\times X$, $Y \subset R^d$ and $X \subset R^{d_x}$.
Assume the conditional distribution $F(y |x)$ has a density $f(y |x)$
with respect to the Lebesgue measure. The marginal density of $x$ with
respect to some generic measure is denoted by $f(x)$.
A model $\mathcal{M}$ for the conditional density $f(y |x)$ is
described by $p(y|x,\mathcal{M})$.
The KL distance between
$f(y|x)f(x)$ and $p(y|x,\mathcal{M})f(x)$ is defined by
\[
d_{\mathrm{KL}}(F,\mathcal{M}) = \int\log\frac{f(y|x)}{p(y|x,\mathcal{M})}
F(dy,dx).
\]
This distance can also be interpreted as the expected KL distance
between the conditional distributions. Either way, this is the distance
useful for obtaining estimation consistency results.
Also, convergence in the KL distance implies convergence in the total
variation distance.
Below, I consider several different specifications of mixture of normal
regressions models, $p(y|x,\mathcal{M})$, and provide conditions on
$F$ under which $d_{\mathrm{KL}}(F,\mathcal{M})$ can be made arbitrarily small.
I also derive rates of convergence and easy to interpret bounds for
$d_{\mathrm{KL}}(F,\mathcal{M})$.

In general, a finite mixture of normal regressions model can be written as
\[
p(y|x,\mathcal{M}) = \sum_{j=1}^m \alpha_j^m(x )
\phi(y, \mu_j^m(x ),\sigma_j^m(x )),
\]
where mixing probabilities satisfy
$
\alpha_j^m(x) \in[0,1]$ and $\sum_j \alpha_j^m(x) = 1$,
and $\phi(y, \mu$, $\sigma)$ is a normal density with mean $\mu$ and
standard deviation $\sigma$ evaluated at $y$ (if $y$ is
multidimensional then the variance--covariance matrix is diagonal
$\sigma
^2 I$).
Most of the results obtained in the paper can be easily extended to
models in which general location scale densities $\sigma^{-d}K((y-\mu
)/\sigma)$ are mixed instead of the normal densities $\phi(y, \mu
,\sigma)$.
Models, in which the mixing weights depend on $x$, are referred in this
paper as smooth mixtures.
In practice, $\alpha_j^m(x)$'s are often modeled by a multinomial choice
model, for example, multinomial logit [\citet
{PengJacobsTanner1996}] or probit
[\citet{Geweke07}], or it might not depend on $x$.
The mean $\mu_j^m(x)$ can be constant, linear or flexible, for example,
polynomial, in $x$. An exponentiated polynomial or spline in $x$ can be
used for modeling the standard deviation $\sigma_j^m(x)$ [\citet
{VillaniKohnGiordani07}].

To the best of my knowledge, previous literature
on smooth mixtures of regressions (or experts) does not provide a
theory on
what specifications for $\alpha_j^m$, $\mu_j^m$ and
$\sigma_j^m$ deliver a model that can approximate and consistently
estimate large nonparametric classes of densities $F$.
There are theoretical results on approximation of smooth functions and
estimation of conditional expectations by ME [see \citet
{ZeeviMeirMaiorov98} and \citet{MaiorovMeir98}].
The only paper on approximation of
conditional densities by ME seems to be \citet{JiangTanner99} who
develop approximation and estimation results for target densities from
a single parameter exponential family, in which the parameter is a
smooth function of covariates. A detailed comparison with results in
\citet{JiangTanner99} is presented in Section \ref{sec:comparison}.
In this paper, I do not restrict the functional form of $f(y|x)$ and
use weak regularity conditions to describe a class of $F$ that can be
approximated.
Conditions on approximable classes of $f(y|x)$ and $f(x)$ that are
common for different
model specifications include
bounded support for $f(x)$,
continuity of $f(y|x)$ in $(y,x)$,
finite expectation of a change of $\log f(y|x)$ in a neighborhood of
$y$ and existence of the second moments of $y$.
The latter restriction can be weakened by adding densities with fat
tails to the mixtures in addition to normal densities.

In Section \ref{sec:lin_logit}, I show that considerable flexibility is
already attained when $\alpha_j^m$'s are modeled by
multinomial logit with linear\vspace*{1pt} indices in
$x$, and $(\mu_j^m,\sigma_j^m)$ are independent of $x$.
Results in Sections \ref{sec:poly_logit} and \ref{sec:lin_logit}
suggest that using polynomials in the logit specification reduces the
number of mixture components $m$ required to achieve a specified
approximation precision. As shown in Section \ref{sec:flex_mean},
models for univariate response $y$ in which
the mixing probabilities and
the variances of the mixed normals are independent of $x$, and the
means are flexible, for example, polynomial in $x$,
can approximate large classes of $f(y|x)$. Differences in quantiles of
$f(y|x)$ from these classes have to be bounded above and below
uniformly in $x$.
These restrictions on $f(y|x)$ can be weakened if
the variances of the mixed normals are modeled by flexible functions of
$x$. Section \ref{sec:conclusion} summarizes the findings.

%s2 ###
\section{Infeasible model}
\label{sec:Infeasible_model}

In this section, I explicitly construct
a smooth mixture of normals model
that converges to a given $F$ in the KL distance as $m$ increases.
This model is not feasible in the sense that it is not based on
components employed in practice, for example, logit/probit mixing
probabilities. However, the results for feasible models presented in
the following sections follow from this one or are similar.

Let $A_j^m$, $j=0,1,\ldots,m$, be a partition of $Y$ consisting of
adjacent half-open half-closed hypercubes $A_1^m,\ldots,A_m^m$ with
side length
$h_m$ and the rest of the space~$A_0^m$.
As $m$ increases the fine part of the partition becomes finer, $h_m
\rightarrow0$.
Also, it covers larger and larger part of $Y$:
for any $y \in Y$
there exists $M_0$ such that
%
%e2.1 ###
%
\begin{equation}
\label{eq:con_partition_gen_case}
\forall m \geq M_0\qquad C_{\delta_m}(y) \cap A_0^m = \varnothing,
\end{equation}
where $C_{\delta_m}(y)$ is a hypercube with center $y$ and side length
$\delta_m \rightarrow0$. It is always possible to construct such a
partition. For example, if $Y=[0,\infty)$ let $A_0^m=[\log m, \infty)$,
$A_j^m=[(j-1)\log m / m, j \log m / m)$ for $j \neq0$, and $h_m=\log
m / m$.

A candidate model $\mathcal{M}_0$ for approximating $f(y|x)$ is
%
%e2.2 ###
%
\begin{equation}
\label{eq:candidate_general_case}
p(y|x,\mathcal{M}_0) = \sum_{j=1}^m F(A_j^m|x) \phi(y, \mu
_j^m,\sigma_m)
+F(A_0^m|x) \phi(y, 0,\sigma_0),
\end{equation}
where $\sigma_0$ is fixed, $\sigma_m$ converges to zero as $m$ increases
and $\mu_j^m$ is the center of $A_j^m$.
One can always construct a model $\mathcal{M}_0$ and a partition
$A_j^m$ so that
%
%e2.3 ###
%
\begin{equation}
\label{eq:cond_delta_sigma_h}
\delta_m \rightarrow0,\qquad \sigma_m / \delta_m \rightarrow0,\qquad
\delta
_m^{d-1} h_m / \sigma_m^d \rightarrow0,
\end{equation}
for example, in the example for $Y=[0,\infty)$ from the previous
paragraph let $\sigma_m = h_m^{0.5}$ and $\delta_m = h_m^{0.25}$.

For a partition satisfying
(\ref{eq:con_partition_gen_case}) and (\ref{eq:cond_delta_sigma_h}),
let us introduce the following restrictions on $F$.
\begin{assumption} \label{assn:general_case}
1.
\hypertarget{assnitem:general_case_1}
$f(y|x)$ is continuous in $y$ a.s. $F$.
%on the interior of the support $Y$ for all $x \in X$.
%
\smallskipamount=0pt
\begin{enumerate}[2.]
\item[2.]
\hypertarget{assnitem:general_case_2}
The second moments of $y$ are finite.
\item[3.]
\hypertarget{assnitem:general_case_3}
For any $(y,x)$ there exists a hypercube $C(r,y,x)$ with side length $r>0$
and $y \in C(r,y,x)$
such that (i)
%
%e2.4 ###
%
\begin{equation}
\label{eq:IntBoundFinfF}
\int\log\frac{f(y|x)}{\inf_{z \in C(r,y,x)} f(z|x) } F(dy,dx) <
\infty
\end{equation}
and (ii) exists $M_3$ such that for any $m \geq M_3$,
if $y \in A_0^m$ then $C(r,y,x) \cap A_0^m$
contains a hypercube $C_0(r,y,x)$ with side length $r/2$ and a vertex
at $y$ and
if $y \in Y \setminus A_0^m$, then $C(r,y,x) \cap(Y \setminus A_0^m)$
contains a hypercube $C_1(r,y,x)$ with side length $r/2$ and a vertex
at $y$.
\end{enumerate}
\end{assumption}

Parameter $\sigma_0$ can always be chosen so that
%
%e2.5 ###
%
\begin{equation}
\label{eq:cond_sigma0}
1>2^{-(d+1)}>\phi(y, 0,\sigma_0) \lambda(C_0(r,y,x)),
\end{equation}
where $\lambda$ is the Lebesgue measure.
\begin{proposition} \label{prp:gen_case}
If the model $p(y|x,\mathcal{M}_0)$ and the partition $A_j^m$ are
constructed so that
(\ref{eq:con_partition_gen_case}), (\ref{eq:candidate_general_case}),
(\ref{eq:cond_delta_sigma_h}) and (\ref{eq:cond_sigma0}) hold, and
$F$ satisfies Assumption \ref{assn:general_case},
then
$d_{\mathrm{KL}}(F,\mathcal{M}_0) \rightarrow0$ as $m\rightarrow\infty$.
\end{proposition}

The proposition is rigorously proved in the \hyperref[app]{Appendix}.
Here, I briefly
describe the intuition behind the argument and the role of the assumptions.
Convergence in the KL distance is proved by the dominated convergence
theorem (DCT). First, I establish point-wise convergence of the integrand,
$\log f(y|x)/ p(y|x ,\mathcal{M}_0)$, to zero, and then I derive an
integrable upper bound on the integrand for the DCT applicability.
Nonnegativity of the KL distance is fruitfully exploited in the proof
as it allows working only with upper bounds and ignoring the lower ones
in convergence arguments.

The first term on the right-hand side of (\ref
{eq:candidate_general_case}) (the sum from 1 to $m$) approximates the integral
%
%e2.6 ###
%
\begin{equation}
\label{eq:convol}
\int\phi(y, \mu,\sigma_m) f(\mu|x) \,d \mu=
\int f(y - \sigma_m z | x)
\phi(z, 0,1) \,d z ,
\end{equation}
when $h_m$ is much smaller than $\sigma_m$, and the fine part of the
partition is large.
The integral on the right-hand side of (\ref{eq:convol}) is obtained by
the change of variables.
For a small $\delta_m$ and $z$ satisfying $\Vert\sigma_m z \Vert\leq
\delta
_m$, $f(y-\sigma_mz|x)$ is close to $f(y|x)$ as $f(y|x)$ is assumed to
be continuous in $y$. Therefore, when $\sigma_m$ is much smaller than
$\delta_m$ the right-hand side of (\ref{eq:convol}) should be close to
$f(y|x)$. Thus, this intuitive argument explains the role of conditions
(\ref{eq:cond_delta_sigma_h}) and continuity of $f(y|x)$.

The second term on the right-hand side of (\ref
{eq:candidate_general_case}) converges to zero. This term is not needed
for point-wise convergence. It can be omitted when the support of
$f(y|x)$ is bounded uniformly in $x$ as in this case we can set
$A_0^m=\varnothing$ and use the same variance $\sigma_m^2$ in all mixture
components (there is no need to define $\sigma_0$).
This term
together with part \hyperlink{assnitem:general_case_2}{2} of Assumption
\ref
{assn:general_case}
prevents tails of $p(y|x ,\mathcal{M}_0)$ from becoming too thin
relative to $f(y|x)$
in the unbounded support case (in the absence of this term the tails
would be too thin as $\sigma_m \rightarrow0$).
%Excessively thin tails of $p(y | x ,\mathcal{M}_0)$ relative to
%$f(y|x)$ are ruled out by part \ref{assnitem:general_case_2} of
%Assumption \ref{assn:general_case}.

Parts \hyperlink{assnitem:general_case_2}{2} and \hyperlink
{assnitem:general_case_3}{3}
of Assumption \ref{assn:general_case} together guarantee existence of
an integrable upper bound for the DCT applicability.
An upper bound on $\log f(y|x) / p(y | x $, $\mathcal{M}_0)$
involves a lower bound on $p(y | x ,\mathcal{M}_0)$. Both terms on the
right-hand side in the definition of $p(y | x ,\mathcal{M}_0)$ in
(\ref{eq:candidate_general_case}) can be bounded below
by an expression proportional to $\inf_{z \in C(r,y,x)} f(z|x)$. That
is how condition (\ref{eq:IntBoundFinfF}) is deduced.
The lower bound for the second term in (\ref
{eq:candidate_general_case}) also includes $\phi(y, 0,\sigma_0)$ and
that is why finiteness of the second
moments of $y$ is assumed.

One interpretation of condition (\ref{eq:IntBoundFinfF})
[part \hyperlink{assnitem:general_case_3}{3}(i) of Assumption \ref
{assn:general_case}]
is that local relative changes in
$f(y|x)$ due to changes in $y$ should not be infinitely large
on average.
%over considerable parts of $Y$ and $X$.
It seems difficult to think of an unconditional density, which is well
behaved and positive everywhere, that would violate (\ref{eq:IntBoundFinfF}).
This part of the assumption though can be violated by reasonable
conditional densities as Example~\ref{ex:exponential} below illustrates.

%
%f1 ###
%
\begin{figure}[b]

\includegraphics{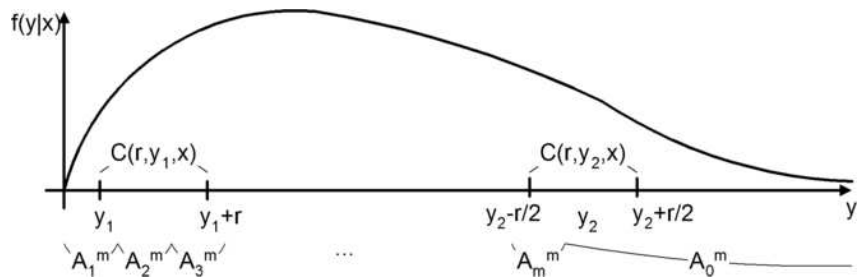}

\caption{Construction of $C(r,y,x)$.}
\label{fig:Cryx}
\end{figure}

When $f(y|x)$ is positive everywhere, part \hyperlink
{assnitem:general_case_3}{3}(ii)
of Assumption \ref{assn:general_case} is not needed.
It always holds if $C(r,y,x)$ is a hypercube with center at $y$.
Part \hyperlink{assnitem:general_case_3}{3}(ii) becomes important when $f(y|x)$
can be equal to zero.
In particular, the sets $C_0(r,y,x)$ and $C_1(r,y,x)$ in
part \hyperlink{assnitem:general_case_3}{3}(ii) of Assumption \ref
{assn:general_case}
are introduced to specify that $C(r,y,x)$ needs to be defined
differently near the boundary of the support and in the tails if one
wants to use condition (\ref{eq:IntBoundFinfF}) in its present form.
This is illustrated in Figure \ref{fig:Cryx}.
%Near the boundary of the support, $y$ has to be at the boundary of
%$C(r,y,x)$ while in the tails, $y$ has to be in the middle of
%$C(r,y,x)$.
%These requirements on construction of $C(r,y,x)$ are employed in
%showing that $\inf_{z \in C(r,y,x)} f(z|x)$ is indeed a part of the
%lower bound on $p(y | x ,\mathcal{M}_0)$ for all sufficiently large
%$m$.
%Figure \ref{fig:Cryx} illustrates
%construction of $C(r,y,x)$ near the boundary of the support and in the
%tails of $p(y|x)$.

The support of $f(\cdot|x)$ should include $C(r,y,x)$ a.s. $F$;
otherwise, part \hyperlink{assnitem:general_case_3}{3}(i)
of Assumption \ref{assn:general_case} is not satisfied.
Therefore, for $f(y|x)$ in Figure \ref{fig:Cryx}, it has to be the case
that $C(r,y,x)=[y,y+r]$ at the boundary of the support (the
intersection of the axes).
Setting $C(r,y,x)=[y,y+r]$ near the boundary of the support
makes the ratio $f(y|x) / \inf_{z \in C(r,y,x)} f(z|x)$ smallest
possible (equal to one) and
thus helps with condition (\ref{eq:IntBoundFinfF}).
Parts of $Y$ near the boundary of the support are covered by the fine
part of the partition $A_1^m,\ldots,A_m^m$ for all sufficiently large
$m$, and
part \hyperlink{assnitem:general_case_3}{3}(ii)
of Assumption \ref{assn:general_case} holds for $C_1(r,y,x)=[y,y+r/2]$.
Using $C(r,y,x)=[y,y+r]$ for all $y$ would not work. Since for any $m$
one can find $y \in A_m^m$ such that $C(r,y,x) \cap Y\setminus A_0^m$
is arbitrary small, and part \hyperlink{assnitem:general_case_3}{3}(ii)
of Assumption \ref{assn:general_case} fails.
Thus, for $y$ that are arbitrary far from the boundary of the support,
one has to use $C(r,y,x)=[y-r/2,y+r/2]$ eventually. Then, part
\hyperlink{assnitem:general_case_3}{3}(ii) of the assumption clearly
holds for
$C_1(r,y,x)=[y-r/2,y]$, $C_0(r,y,x)=[y,y+r/2]$ and any $m$.

%Part \ref{assnitem:general_case_3}(ii)
%of Assumption \ref{assn:general_case} also rules out cases when the
%support of $f(\cdot|x)$ increases without a bound as $x$ changes,
%e.g., $X=[0,\infty)$ and the support of $f(\cdot|x)$ is $[-x,x]$.
%Cases in which the support boundary is never completely covered by the
%fine part of the partition are not necessarily ruled out, e.g.,
%an example with a two dimensional $y$ with support of $f(\cdot|x)$
%equal to
%$[0,\infty)\times[0,\infty)$ can be accommodated.

Results in this section and similar results in the following sections
can be generalized in several different ways.
First, the derivation of the integrable upper bound in the proof of
Proposition \ref{prp:gen_case} suggests that the requirement of finite
second moments of $y$ can be weakened by adding a density with thicker
than normal tails to the mixture of normals; for example, substitute
$\phi(y, 0,\sigma_0)$ in (\ref{eq:candidate_general_case}) with a
Student $t$-density.
Second, more general shapes of the support of $F$
can be accommodated if
instead of hypercubes $C(r,y,x)$, $C_0(r,y,x)$, and $C_1(r,y,x)$ in
Assumption \ref{assn:general_case}
different sets with positive Lebesgue measure are used.
For example, if the support of $f(\cdot|x)$ is a triangle in $R^2$ then
small triangles can be used instead of the squares $C(r,y,x)$,
$C_0(r,y,x)$ and $C_1(r,y,x)$.
%It seems possible to weaken part \ref{assnitem:general_case_3}(ii)
%of Assumption \ref{assn:general_case} to accommodate some special
%cases in which
%the support boundary of $f(\cdot|x)$ is not bounded uniformly in $x$.
%However, this generalization is not pursued in the paper since
%model specifications in which normal means and variances are functions
%of $x$ seems to be more suitable as they can accommodate more general
%situations; these specifications are considered in Section
Third, general location scale densities $\sigma^{-d}K((y-\mu)/\sigma)$
can be used in mixtures instead of normal densities.
As long as analogs of Lemmas \ref{lm:boundRSbyInt_ESP}, \ref
{lm:boundIntBy1} and \ref{lm:boundRS_EPP_new} (see the \hyperref
[app]{Appendix}) are
available for a particular type of densities,
results in this and the following sections
will hold for mixtures of these densities.
Lemmas \ref{lm:boundRSbyInt_ESP} and \ref{lm:boundRS_EPP_new} hold for
$\sigma^{-d}K((y-\mu)/\sigma)$ if
$K(z)$ is bounded and nonincreasing in $|z|$ (proofs of the lemmas use
only these facts about the normal distributions).
The derivation of bounds in Lemma \ref{lm:boundIntBy1} exploits
normality; however, the qualitative results of the
lemma hold as long as $\int_R K(z)\,dz = 1$ and $K(z)$ is positive in a
neighborhood of zero. Thus, all the results in this paper that establish
$d_{\mathrm{KL}}(F,\mathcal{M}) \rightarrow0$
do not depend on the normality assumption; however, bounds and
convergence rates for $d_{\mathrm{KL}}(F,\mathcal{M})$ derived below are
specific to mixtures of normal densities, and they might be different for
mixtures of other densities.
All these generalizations seem to be straight forward
and I do not pursue them in this paper
to keep the arguments short and simple.

Examples below demonstrate that Assumption \ref{assn:general_case} is
satisfied for a large class of densities. They also describe some
situations in which the assumption fails.
%
%Densities that satisfy Assumption \ref{assn:general_case}.
%
\begin{example}
\label{ex:exponential}
Exponential distribution, $f(y|x)=\gamma(x) \exp\{-\gamma(x)y\}$,\break
\mbox{$\gamma(x)>0$}.
The density is continuous in $y$ (part \hyperlink
{assnitem:general_case_1}{1} of
Assumption \ref{assn:general_case}).
Let $\int\gamma^{-2} \,dF < \infty$ so that the second moment of $y$ is
finite (part \hyperlink{assnitem:general_case_2}{2} of Assumption~\ref
{assn:general_case}).
Define the partition $A_j^m$ and $C(r,y,x)$, $C_0(r,y,x)$ and
$C_1(r,y,x)$ as shown in Figure \ref{fig:Cryx}, for example,
for some $r>0$ let $C(r,y,x)=[y,y+r]$ for $y \in[0,r]$ and
$C(r,y,x)=[y-r/2,y+r/2]$ for $y \in(r,\infty)$.
Thus, from the discussion of Figure \ref{fig:Cryx} above it follows
that part \hyperlink{assnitem:general_case_3}{3}(ii) of Assumption \ref
{assn:general_case} is satisfied.
Because $\log f(y|x) / \inf_{z \in C(r,y,x)} f(z|x) \leq r \gamma(x)$,
part \hyperlink{assnitem:general_case_3}{3}(i) of Assumption \ref
{assn:general_case} holds as long as $\gamma(x)$ is integrable with
respect to $f(x)$. If $\gamma(x)$ is not integrable, then part
\hyperlink{assnitem:general_case_3}{3}(i) of the assumption fails.
\end{example}
\begin{example}
\label{ex:inf_student}
A Student $t$-distribution, in which scale and location parameters are
functions of $x$,
$f(y|x) \propto[ \nu+ ((y-b(x))/c(x))^2]^{-(\nu+ 1)/2}$, $\nu>2$ and
$b(x)^2$, $c(x)^{-2}$ and $c(x)^2$ are integrable w.r.t. $f(x)$.
The second moment of $y$ is finite since
\begin{eqnarray*}
\int y^2 \,dF & = & \int\biggl(c(x)^2 \biggl[ \frac{y-b(x)}{c(x)} \biggr]^2 + 2
b(x)y - b(x)^2 \biggr) \,dF \\
& = & \int\biggl( c(x)^2 \frac{\nu}{\nu- 2} + 2 b(x)^2 - b(x)^2
\biggr) \,dF < \infty.
\end{eqnarray*}
As I discuss above, for densities positive everywhere
part \hyperlink{assnitem:general_case_3}{3}(ii) of Assumption \ref
{assn:general_case} always holds with $C(r,y,x) = [y-r/2,y+r/2]$.
Part \hyperlink{assnitem:general_case_3}{3}(i) of Assumption \ref
{assn:general_case} is also satisfied because
\begin{eqnarray*}
&& \int\log\frac{f(y|x)}{\inf_{z \in C(r,y,x)} f(z|x) } F(dy,dx) \\
&&\qquad = 2 \int_X \int_{b(x)}^{\infty} -\frac{\nu+ 1}{2}
\log\frac{\nu+ ((y-b(x))/c(x))^2}{\nu+ ((y+ r - b(x))/c(x) )^2}
f(y|x)\, dy F(dx) \\
&&\qquad \leq(\nu+ 1) 2 \int_X \int_{b(x)}^{\infty}
\bigl[\nu+ \bigl(\bigl(y + r - b(x)\bigr)/c(x) \bigr)^2\bigr] f(y|x)\, dy F(dx) <
\infty,
\end{eqnarray*}
where the last inequality follows by the integrability of
$((y - b(x))/c(x)$, its square and $c(x)^{-2}$.
\end{example}
\begin{example}
\label{ex:inf_cont_bddsupport}
Suppose that conditional density $f(y|x)$ is continuous in $y$ and
bounded above and away from zero,
$\infty> \overline{f} \geq f(y|x) \geq\underline{f}>0$ for any $y
\in
Y=[a,b]$ and $x \in X$.
Then we can set $A_0^m=\varnothing$. For $r \in(0, (b-a)/4)$, let
$C(r,y,x) = [y, y + r]$ and $C_1(r,y,x) = [y, y + r/2]$ for $y \in
[a,(a+b)/2]$ and
$C(r,y,x) = [y-r, y]$ and $C_1(r,y,x) = [y-r/2, y]$ for $y \in((a+b)/2,b]$.
Clearly, part \hyperlink{assnitem:general_case_3}{3}(ii) of Assumption
\ref
{assn:general_case} is satisfied.
Because $f(y|x)/ \inf_{z \in C(r,y,x)} f(z|\break x) \leq\overline{f} /
\underline{f}$ part \hyperlink{assnitem:general_case_3}{3}(i) of Assumption
\ref
{assn:general_case} also holds.
The second moment of $y$ is finite and thus all parts of Assumption
\ref
{assn:general_case} hold.

The boundedness away from zero condition can be replaced by a monotonicity condition
at the boundary of the support. For example, let $f(y|x)$ be nondecreasing on $[a, a + 2r]$,
nonincreasing on $[b - 2r, b]$ and bounded below by $\underline{f} > 0$ on $[a + r, b - r]$. In this case
$f(y|x)/ \operatorname{inf}_{z\in C(r,y,x)} f(z|x) \leq \operatorname{max}\{1, \overline{f}/\underline{f}\}$ for any $y \in [a, b]$. Thus, part 3(i) of Assumption
\ref{assn:general_case} holds. The other parts of the assumption are not affected by this change.

%Bounded away from zero condition can be substituted with monotonicity
%at the boundary of the support; for example,
%$f(y|x)$ is nondecreasing on $[a,a+2 r]$, nonincreasing on $[b-2r,b]$
%and bounded below by $\underline{f}>0$ on $[a+r,b-r]$,
%so that
%$f(y|x)/ \inf_{z \in C(r,y,x)} f(z|x) \leq\max\{1,\overline{f} /
%holds. The other parts of the assumption are not affected by this change.
\end{example}
\begin{example}
\label{ex:infeasible_unifrm}
Consider a uniform distribution $f(y|x) = x^{-1} 1_{[0,x]}(y)$ and
$f(x)>0$ for any $x \in[1,\infty)$.
A natural choice of the partition would be $A_0^m=[m h_m, \infty)$ and
$A_j^m=[(j-1)h_m,j h_m)$ for $j \in\{1,\ldots,m\}$.
When $y=x$, the only reasonable choice of $C(r,y,x)$ is
$C(r,y,x)=[y-r,y]$. For an arbitrary $m$ and $y=x=m h_m+r/4$,
$C(r,y,x)$ violates part \hyperlink{assnitem:general_case_3}{3}(ii) of
Assumption \ref{assn:general_case} since the only possible $C_0(r,y,x)
= [y-r/2,y]$ is not included in $A_0^m$.
For $f(x)$ with bounded support, this example would satisfy Assumption
\ref{assn:general_case} since in this case we could set
$A_0^m=\varnothing$.

This example illustrates that Assumption \ref{assn:general_case} rules
out some cases in which the support of $f(\cdot|x)$ is increasing in
$x$ without a bound. In Section \ref{sec:flex_mean}, I consider model
specifications in which means and variances of the mixed normals can be
flexible functions of $x$. Those specifications seem to be more
promising for modeling densities $f(\cdot|x)$ with support increasing
in $x$ without a bound (see Example \ref{ex:uniform_flex_mean}).
\end{example}

%s2.1 ###
\subsection{Approximation error bounds}
\label{sec:M0bounds}

The proof techniques of this section can also be used to derive
explicit bounds on the approximation error.
The bounds for positive everywhere and especially differentiable
$f(y|x)$ are particularly informative.
It is also easy to deduce an approximation rate from them.
Thus, I present below the bounds and approximation rate for these
special albeit important cases.
Convergence rates and bounds for other special classes can be obtained
in a similar way, for example, for densities bounded away from zero.
However, rates and bounds for the general case seem to be difficult to
calculate.
\begin{corollary} \label{crl:gen_case_bounds}
Part \textup{(i)}. Suppose
the model $p(y|x ,\mathcal{M}_0)$ and the partition $A_j^m$ are
constructed so that
(\ref{eq:con_partition_gen_case}), (\ref{eq:candidate_general_case}),
(\ref{eq:cond_delta_sigma_h}) and\vspace*{-2pt} (\ref{eq:cond_sigma0}) hold.
Suppose $f(y|x)$ is positive and continuous in $y$ on $Y=R^d$ for all
$x$, second moments of $y$ are finite and
(\ref{eq:IntBoundFinfF}) holds with $C(r,y,x)=C_r(y)$ taken to be a
hypercube with center at $y$ and radius $r$.
Then, for all sufficiently large $m$,
%
%e2.10 ###
%e2.9 ###
%e2.8 ###
%e2.7 ###
%
\begin{eqnarray}
\label{eq:bound1_intFinfF}
d_{\mathrm{KL}}(F,\mathcal{M}_0)
&\leq& \int\log\frac{f(y|x)}{\inf_{z \in C_{\delta_m}(y)} f(z|x)} F(dy,dx)
\\
\label{eq:bound2_lof1_epsm}
&&{} + 2 \frac{3 d^{3/2} \delta_m^{d-1} h_m }{(2 \pi)^{d/2} \sigma_m^d}
+ 2 \exp\biggl\{-\frac{(\delta_m/\sigma_m)^2}{8} \biggr\}
\\
\label{eq:bound3_intFinfFtail}
&&{} + \int_{B_{\delta_m}(A_0^m)}
\log\frac{f(y|x)}{\inf_{z \in C_{r}(y)} f(z|x)} F(dy,dx)
\\
\label{eq:bound4_inty2_cFtail}
&&{} + \int_{B_{\delta_m}(A_0^m)}
\biggl[ \frac{y^{\prime}y}{2 \sigma_0^2} -
\log\frac{(r/2)^d}{(2 \pi\sigma_0^2)^{d/2}} \biggr] F(dy,dx),
\end{eqnarray}
where
$B_{\delta_m}(A_0^m)=\{(y,x)\dvtx
C_{\delta_m}(y)\cap A_0^m \neq\varnothing
\}
$%\]
and bounds in (\ref{eq:bound1_intFinfF})--(\ref{eq:bound4_inty2_cFtail})
converge to zero as $m \rightarrow\infty$.

Part \textup{(ii)}. If $f(y|x)$ is continuously differentiable in $y$ for all
$x$ and instead of (\ref{eq:IntBoundFinfF}) the following condition holds:
%
%e2.11 ###
%
\begin{equation}
\label{eq:dfdy_integrable}
\int\sup_{z \in C_{r}(y)} \biggl\Vert\frac{d \log f(z|x) }{dz}\biggr\Vert F(dy,dx) <
\infty,
\end{equation}
then for all sufficiently large $m$,
%
%e2.15 ###
%e2.14 ###
%e2.13 ###
%e2.12 ###
%
\begin{eqnarray}\qquad
\label{eq:bound1_intdsuplnFdz}
d_{\mathrm{KL}}(F,\mathcal{M}_0)
&\leq& \delta_m \cdot\frac{d^{1/2}}{2} \int\sup_{z \in C_{\delta
_m}(y)} \biggl\Vert\frac{d \log f(z|x) }{dz}\biggr\Vert F(dy,dx)
\\
\label{eq:bound2_lof1_epsm2}
&&{} +2 \frac{3 d^{3/2} \delta_m^{d-1} h_m }{(2 \pi)^{d/2} \sigma_m^d}
+2 \exp\biggl\{-\frac{(\delta_m/\sigma_m)^2}{8} \biggr\}
\\
\label{eq:bound3_intdsuplnFdztail}
&&{} + \frac{r d^{1/2}}{2} \int_{B_{\delta_m}(A_0^m)}
\sup_{z \in C_{r}(y)} \biggl\Vert\frac{d \log f(z|x) }{dz}\biggr\Vert F(dy,dx)
\\
\label{eq:bound4_inty2_cFtail_2}
&&{} + \int_{B_{\delta_m}(A_0^m)}
\biggl[ \frac{y^{\prime}y}{2 \sigma_0^2} -
\log\frac{(r/2)^d}{(2 \pi\sigma_0^2)^{d/2}} \biggr] F(dy,dx),
\end{eqnarray}
and bounds in (\ref{eq:bound1_intdsuplnFdz})--(\ref
{eq:bound4_inty2_cFtail_2}) converge to zero as $m \rightarrow\infty$.

Part \textup{(iii)}.
If, in addition to assumptions from part \textup{(ii)}, for some $q>2$ and some
$i_1 \in\{1,\ldots,d\}$
%
%e2.16 ###
%
\begin{equation}
\label{eq:bdd_q_moment}
\int|y_i|^q F(dy) < \infty,\qquad i \in\{1,\ldots,d\},
\end{equation}
and
%
%e2.17 ###
%
\begin{equation}
\label{eq:bdd_y_dlogfdz}
\int|y_{i_1}|^{q-2} \sup_{z \in C_{r}(y)} \biggl\Vert\frac{d \log f(z|x)
}{dz}\biggr\Vert F(dy,dx) < \infty,
\end{equation}
then the approximation error bound can be written as
%
%e2.18 ###
%
\begin{equation}
\label{eq:bd_rate}
d_{\mathrm{KL}}(F,\mathcal{M}_0)
\leq
c \cdot\biggl(\frac{1}{m} \biggr)^{1/(d \cdot[ 2+ 1 / (q-2) +
\varepsilon])},
\end{equation}
where $\varepsilon>0$ can be arbitrarily close to zero and $c$ does not
depend on $m$.
\end{corollary}

The corollary is proved in the \hyperref[app]{Appendix}. The bounds in
part (i) of the
corollary follow from the proof of Proposition \ref{prp:gen_case}.
The bounds in part (ii) are derived from the bounds in part (i),
and they are especially easy to interpret.
The larger the ``average'' derivative of $\log f(y|x)$ is the smaller
$\delta_m$ has to be to achieve a prespecified level for the right-hand
side of (\ref{eq:bound1_intdsuplnFdz}).
Constant $h_m$ has to be much smaller than $\sigma_m$, and $\sigma_m$
has to be much smaller than $\delta_m$ [condition (\ref
{eq:cond_delta_sigma_h})]
so that (\ref{eq:bound2_lof1_epsm2}) becomes sufficiently small. Size
of (\ref{eq:bound3_intdsuplnFdztail}) and (\ref
{eq:bound4_inty2_cFtail_2}) depends on how fast and by how much
tails of $f(y|x)f(x)$ dominate $d \log f(y|x)/dy$, $y^2$, and a constant.

The approximation rate in part (iii) is derived from the bounds in part (ii).
Expressions in (\ref{eq:bound1_intdsuplnFdz}) and (\ref{eq:bound2_lof1_epsm2})
can be immediately converted in expressions in terms of $m$.
To convert
(\ref{eq:bound3_intdsuplnFdztail}) and
(\ref{eq:bound4_inty2_cFtail_2}) in expressions in terms of $m$ one
seems to need slightly more than
integrability of ${\sup_{z \in C_{r}(y)}} \Vert d \log f(z|x) / dz\Vert$
[condition (\ref{eq:bdd_y_dlogfdz})] and slightly more than finiteness
of the second moments of $y$ [condition
(\ref{eq:bdd_q_moment})].
Under these conditions, (\ref{eq:bound3_intdsuplnFdztail}) and
(\ref{eq:bound4_inty2_cFtail_2}) are bounded by $(h_m
m^{1/d})^{-(q-2)}$ times a constant (see the corollary proof).
An upper bound on $(h_m m^{1/d})^{-(q-2)}$, (\ref
{eq:bound1_intdsuplnFdz}) and (\ref{eq:bound2_lof1_epsm2}) gives
the rate in (\ref{eq:bd_rate}). This upper bound
has to be strictly larger than
(\ref{eq:bd_rate}) with $\varepsilon=0$ as I show
in the corollary proof.
For distributions with exponentially declining tails, (\ref
{eq:bound3_intdsuplnFdztail}) and
(\ref{eq:bound4_inty2_cFtail_2}) can be decreasing exponentially in
$h_m m^{1/d}$. In this case, one can set $q=\infty$ in (\ref
{eq:bd_rate}) (see Example \ref{ex:2side_exp_flex_mean_rates} below).

The dimension of $y$ enters the approximation bounds exponentially. The
dimension of $x$ does not affect the bound and the approximation rate
for the ``infeasible'' model because this model is constructed with the
use of $F(A_j^m|x)$'s, which are unknown functions of $x$. The following
sections shed some light on the role of the dimension of $x$ in
approximating $f(y|x)$ by feasible models.

%s3 ###
\section{Flexible multinomial choice models for mixing probabilities}
\label{sec:poly_logit}

This section gives conditions under which
approximation results for
``infeasible'' model $\mathcal{M}_0$ also hold for
a model with logit mixing probabilities that include polynomial terms
in $x$. It also shows how to extend these results to multinomial probit
and other models for mixing probabilities.
\begin{assumption} \label{assn:XcompactLogFcont}
$X$ is compact and for partitions $A_j^m$, $j=0,1,\ldots,m$ satisfying
(\ref{eq:con_partition_gen_case}), $F(A_j^m|x)$ is a continuous
function of $x$ on $X$ and $F(A_j^m|x)>0$ [the support of $f(\cdot|x)$
does not depend on $x$].
\end{assumption}

Under this assumption (by the Stone--Weierstrass theorem) for any
sequence of $\varepsilon_m \rightarrow0$, $\varepsilon_m>0$ there exist
finite order polynomials in $x$, $P_j^m(x )$ such that
%
%e3.1 ###
%
\begin{equation}
\label{eq:logFcont}
|{\log F(A_j^m|x) - P_j^m(x )}| < \varepsilon_m\qquad \forall x \in X,
j=1,\ldots,m.
\end{equation}
Let $p(y|x, \mathcal{M}_1)$ denote a model with $\sigma_j^m$ and $\mu
_j^m$ independent of $x$ and logit mixing probabilities,
\begin{eqnarray*}
\alpha_j^m(x ,\mathcal{M}_1) & = & \frac{ \exp\{ P_j^m(x ) \}}{\sum
_{k=1}^m \exp\{ P_k^m(x )\} } \nonumber\\
& = & \frac{ F(A_j^m|x) \exp\{ P_j^m(x ) - \log F(A_j^m|x) \}}{\sum
_{k=1}^m F(A_k^m|x) \exp\{ P_k^m(x ) - \log F(A_k^m|x)\} }.
\end{eqnarray*}
Condition (\ref{eq:logFcont}) implies
$\alpha_j^m(x,\mathcal{M}_1) \in(F(A_j^m|x) \exp\{-2 \varepsilon_m\},
F(A_j^m|x) \exp\{2 \varepsilon_m\})$. The following corollary
immediately follows.
\begin{corollary} \label{crl:logit_polynomials}
If Assumption \ref{assn:XcompactLogFcont} and the conditions of
Proposition \ref{prp:gen_case} hold then
$d_{\mathrm{KL}}(F,\mathcal{M}_1)$
is bounded above and below by $d_{\mathrm{KL}}(F,\mathcal{M}_0) \pm2 \varepsilon
_m$ and thus converges to zero.
\end{corollary}

It seems possible to extend this corollary to other models
for mixing probabilities, in particular, to a class of
multinomial choice models
in which mixing probabilities have the following representation:
\[
\alpha_j^m(x)=\operatorname{Pr} [(e_0,\ldots,e_m)\dvtx v_j(x)+e_j \geq v_k(x)+e_k,
k \in\{0,\ldots,m\} ],
\]
where $v_j(x)$ are flexible functions of $x$ and
$e_k$'s are i.i.d.
Multinomial logit and probit models fall into this category with
polynomial $v_j(x)$ and extreme value and normal distributions for $e_k$'s.
The proof of Proposition 1 in
\citet{HotzMiller93} implies that if
$e_k$ are i.i.d. and have a density with respect to the Lebesgue
measure, which is positive on $R$, then
\[
(v_0(x),\ldots,v_{m-1}(x)) = Q(\alpha_0^m(x),\ldots,\alpha_{m-1}^m(x)),
\]
where $v_{m}(x)$ is normalized to 0 and $Q$ and $Q^{-1}$ are
differentiable mappings defined correspondingly on $R^m$, and the
interior of the $m$-dimensional simplex.
Flexible functional forms for $(v_0(x),\ldots,v_{m-1}(x))$ can be used
to approximate
$Q(F(A_0^m|x),\ldots,F(A_{m-1}^m|x))$. Then
$(\alpha_0^m(x),\ldots,\alpha_{m-1}^m(x)) =\break Q^{-1}(v_0(x),\ldots
,$ $v_{m-1}(x))$ will approximate
$(F(A_0^m|x),\ldots,F(A_0^{m-1}|x)$.
To get an analog of Corollary \ref{crl:logit_polynomials} one only
needs to show that $Q^{-1}$
transfers
small additive approximation errors in $v_j(x)$ into
multiplicative\vspace*{-1pt} approximation errors for $\alpha_j^m(x)$, that are close
to one.
Since the mapping $Q^{-1}$ is continuous this is the case as long as
$F(A_j^m|x)$ are positive.
Thus, it seems
one does not need more than
Assumption \ref{assn:XcompactLogFcont}
to extend Corollary \ref{crl:logit_polynomials} to other models for
mixing probabilities.

Of course, Corollary \ref{crl:logit_polynomials} can be formulated for
any other method for approximating continuous functions in the sup norm
on compacts, for example, for splines instead of the polynomials in the
logit mixing probabilities.

The corollary implies that for $F$ satisfying conditions of Corollary
\ref{crl:gen_case_bounds}, bounds on the approximation error for model
$\mathcal{M}_1$ are given by the bounds in the corollary for $\mathcal
{M}_0$ plus $\varepsilon_m$.
Results from the function approximation theory [see, e.g.,
Section~3.3 in \citet{Rustndp96} for a survey] suggest
that
to achieve a worst case approximation bound $\varepsilon_m$,
computable approximations to Lipschitz continuous functions must
involve the number of parameters proportional to $\varepsilon_m^{-d_x}$
($\varepsilon_m^{-d_x/n}$ if the function has bounded derivatives up to
order $n+1$).
Thus, the number of parameters in the polynomials (or splines) $P_j^m(x)$
depends at best exponentially on the dimension of~$x$.

It might be very difficult to estimate a model with high order
polynomials in the logit mixing probabilities. The following section
shows that it is not necessary to use high order polynomials in logit
specification to attain flexibility. However, as I discuss at the end
of the following section, polynomial terms might reduce the number of
mixture components required to achieve a specified approximation precision.

%s4 ###
\section{Linear indices in logit}
\label{sec:lin_logit}

In this section I explore an alternative approximation to $F(A_j^m|x)$
based on logit mixing probabilities that use only linear indices in~$x$.
The following assumption is a slightly stricter analog of Assumption
\ref{assn:general_case}.
\begin{assumption} \label{assn:logitlin_case}
1. $X=[0,1]^{d_x}$ (the arguments would go through for a bounded $X$).
\smallskipamount=0pt
\begin{enumerate}[2.]
\item[2.]
\hypertarget{assnitem:logitlin_case_1}
$f(y|x)$ is continuous in $(y,x)$ a.s. $F$.%on $Y \times X$.
\item[3.]
\hypertarget{assnitem:logitlin_case_2}
The second moments of $y$ are finite.
\item[4.]
\hypertarget{assnitem:logitlin_case_3}
For any $(y,x)$ there exists a hypercube $C(r,y,x)$ with side length $r>0$
and $y \in C(r,y,x)$
such that (i)
%
%e4.1 ###
%
\begin{equation}
\label{eq:IntBoundFinfFyx}
\int\log\frac{f(y|x)}{\inf_{z \in C(r,y,x), \Vert t-x\Vert\leq r}
f(z|t) } F(dy,dx) < \infty
\end{equation}
and (ii) exists $M$ such that for any $m \geq M$,
if $y \in A_0^m$ then $C(r,y,x) \cap A_0^m$
contains a hypercube $C_0(r,y,x)$ with side length $r/2$ and a vertex
at $y$ and
if $y \in Y \setminus A_0^m$, then $C(r,y,x) \cap(Y \setminus A_0^m)$
contains a hypercube $C_1(r,y,x)$ with side $r/2$ and a vertex at $y$.
\end{enumerate}
\end{assumption}

Let $B_i^m$, $i=1,\ldots,N(m)$ be equal size half-open half-closed
hypercubes forming a partition of $X=[0,1]^{d_x}$. The partition
becomes finer as $m$ increases, $\lambda(B_i^m)=N(m)^{-1} \rightarrow
0$. Let $x_i^m$ denote the center of $B_i^m$.
Before looking at logit let us consider an ``infeasible'' model
$\mathcal{M}_2$,
\[
p(y|x ,\mathcal{M}_2) = \sum_{i=1}^{N(m)}
\Biggl[ \sum_{j=1}^m
\alpha_{ij}^m(x ,\mathcal{M}_2)\phi(y, \mu_j^m,\sigma_m)
+\alpha_{i0}^m(x ,\mathcal{M}_2)\phi(y, 0,\sigma_0) \Biggr],
\]
where the mixing probabilities
$\alpha_{ij}^m(x ,\mathcal{M}_2) = 1_{B_i^m}(x) F(A_j^m|x_i^m)$.
As the partition of $X$ becomes finer,\vspace*{-2pt} model $\mathcal{M}_2$
approximates $\mathcal{M}_0$ because
$F(A_j^m|x) \approx\sum_{i=1}^{N(m)} 1_{B_i^m}(x) F(A_j^m|x_i^m)$
under continuity\vspace*{1pt} of $f(y|x)$ in $x$ (part \hyperlink
{assnitem:logitlin_case_1}{2}
of Assumption~\ref{assn:logitlin_case}).
Since, $\mathcal{M}_2$ is not interesting on its own I do not make this
argument precise here. Instead I employ this idea to get approximation
results for model $\mathcal{M}_3$ constructed similarly to $\mathcal
{M}_2$ but with logit mixing probabilities,
%
%e4.2 ###
%
\begin{eqnarray}\label{eq:appr1Bxi}
\alpha_{ij}^m(x ,\mathcal{M}_3) & = & \frac{
\exp\{ \log F(A_j^m|x_i^m) - R_m (x_i^{m \prime} x_i^m - 2 x_i^{m
\prime
} x ) \}}
{\sum_{k,l} \exp\{\log F(A_k^m|x_l^m) - R_m (x_l^{m \prime} x_l^m - 2
x_l^{m \prime} x )\} } \nonumber\\[-8pt]\\[-8pt]
& = & F(A_j^m|x_i^m) \frac{
\exp\{ - R_m (x_i^{m \prime} x_i^m - 2 x_i^{m \prime} x ) \}}
{\sum_{l} \exp\{ - R_m (x_l^{m \prime} x_l^m - 2 x_l^{m \prime} x
)\} }.\nonumber
\end{eqnarray}
In this expression, $R_m$ is a positive diverging to infinity sequence
that satisfies the following condition:
%
%e4.3 ###
%
\begin{equation}
\label{eq:CondSmRm}
\exp\{-R_m s_m\}/s_m^{d_x/2} \rightarrow0\qquad \mbox{where } s_m = d_x
\lambda(B_i^m)^{2/d_x} \rightarrow0,
\end{equation}
is the squared diagonal of $B_i^m$. This condition
specifies that $R_m$ should increase fast relative to how fine the
partition of $X$ becomes.
It is always possible to define sequence $R_m$ satisfying (\ref
{eq:CondSmRm}), for example, $R_m=s_m^{-2}$.
\begin{proposition} \label{prp:gen_linear_logit}
If condition (\ref{eq:CondSmRm}), Assumption \ref{assn:logitlin_case},
and conditions of Proposition \ref{prp:gen_case} hold then
$d_{\mathrm{KL}}(F,\mathcal{M}_3) \rightarrow0$ as $m\rightarrow\infty$.
\end{proposition}

The proposition is proved in the \hyperref[app]{Appendix}. The proof
shows that
the expression in (\ref{eq:appr1Bxi})
multiplying $F(A_j^m|x_i^m)$ behaves like $1_{B_{i}^m}(x)$ when $R_m$
becomes large and then uses the same arguments as in the proof of
Proposition \ref{prp:gen_case}.
Attempts to develop similar results for mixing probabilities modeled by
multinomial probit [see, e.g., \citet{Geweke07} for applications]
were not successful.
It would not be hard to make multinomial probit mixing probabilities
behave like indicator functions. However, making them behave like an
indicator times $F(A_j^m|x_i^m)$ as in (\ref{eq:appr1Bxi}) seems to be
more difficult.

The bounds on the approximation error for
$\mathcal{M}_3$ and $f(y|x)$ positive everywhere are
similar to bounds for $\mathcal{M}_0$ obtained in Corollary \ref
{crl:gen_case_bounds}.
This is formalized in the following corollary.
\begin{corollary} \label{crl:linlogit_case_bounds}
Part \textup{(i)}. Suppose conditions of Proposition \ref{prp:gen_linear_logit} hold,
$f(y|x)$ is positive for any $y \in Y=R^d$ and any $x \in X$,
$f(y|x)$ is continuously differentiable in $(y,x)$, and instead of
(\ref{eq:IntBoundFinfFyx}) the following condition holds:
%
%e4.4 ###
%
\begin{equation}
\label{eq:dfdyx_integrable}
\int\sup_{y \in C_{r}(y), \Vert x-t\Vert\leq r} \biggl\Vert\frac{d \log f(z|t)
}{d(z,t)}\biggr\Vert F(dy,dx) < \infty;
\end{equation}
then, for all sufficiently large $m$,
%
%e4.9 ###
%e4.8 ###
%e4.7 ###
%e4.6 ###
%e4.5 ###
%
\begin{eqnarray}\label{eq:bound1_intdsuplnFdyx}\qquad\quad
d_{\mathrm{KL}}(F,\mathcal{M}_3)
&\leq&
\biggl(\delta_m \frac{d^{1/2}}{2} +
s_m^{1/2}\biggr)\nonumber\\[-8pt]\\[-8pt]
&&{} \times\int\sup_{z \in
C_{\delta_m}(y), \Vert x-t\Vert\leq s_m^{1/2}} \biggl\Vert\frac{d \log f(z|t)
}{d(z,t)}\biggr\Vert F(dy,dx)
\nonumber\\
\label{eq:bound2_lof1_epsm2_2}
&&{}
+2 \frac{3 d^{3/2} \delta_m^{d-1} h_m }{(2 \pi)^{d/2} \sigma_m^d}
+2 \exp\biggl\{-\frac{(\delta_m/\sigma_m)^2}{8} \biggr\}
\\
\label{eq:bound3_intdsuplnFdzxtail}
&&{} + \frac{r d^{1/2}}{2} \int_{B_{\delta_m}(A_0^m)}
\sup_{z \in C_{r}(y), \Vert x-t\Vert\leq r} \biggl\Vert\frac{d \log f(z|t)
}{d(z,t)}\biggr\Vert F(dy,dx)
\\
\label{eq:bound4_inty2_cFtail_2_2}
&&{} + \int_{B_{\delta_m}(A_0^m)}
\biggl[ \frac{y^{\prime}y}{2 \sigma_0^2} -
\log\frac{(r/2)^d}{(2 \pi\sigma_0^2)^{d/2}} \biggr] F(dy,dx)
\\
\label{eq:bound5_expRmsm}
&&{} + \log[1-d_x^{d_x/2}\exp\{-R_m s_m\}/s_m^{d_x/2} ],
\end{eqnarray}
and bounds in (\ref{eq:bound1_intdsuplnFdyx})--(\ref
{eq:bound5_expRmsm}) converge to zero as $m \rightarrow\infty$.

Part \textup{(ii)}.
If, in addition to assumptions from part \textup{(i)}, for some $q>2$ and some
$i_1 \in\{1,\ldots,d\}$,
%
%e4.10 ###
%
\begin{equation}
\label{eq:bdd_q_moment_logit}
\int|y_i|^q F(dy) < \infty,\qquad i \in\{1,\ldots,d\},
\end{equation}
and
%
%e4.11 ###
%
\begin{equation}
\label{eq:bdd_y_dlogfdzt}
\int|y_{i_1}|^{q-2} \sup_{z \in C_{r}(y), \Vert x-t\Vert\leq r} \biggl\Vert
\frac{d
\log f(z|t) }{d(z,t)}\biggr\Vert F(dy,dx) < \infty,
\end{equation}
then the approximation error bound can be written as
%
%e4.12 ###
%
\begin{equation}
\label{eq:bd_rate_linlogit}
d_{\mathrm{KL}}(F,\mathcal{M}_3) \leq
\mbox{constant} \cdot[m N(m) ]^{-1/ (d_x+d \cdot[2+1 /
(q-2)+\varepsilon])},
\end{equation}
where $m N(m)+1$ is the number of mixture components in
$\mathcal{M}_3$ and $\varepsilon>0$ can be arbitrarily close to zero.
\end{corollary}

%In particular,
%the bounds for $\mathcal{M}_3$ include
%They also include
%slightly modified \eqref{eq:bound1_intFinfF}
%and
%$y$ but also over $x$ (compare part 3.(i) of Assumptions
%For differentiable $f(y|x)$, the derivatives and sup's in
%and \eqref{eq:bound3_intdsuplnFdztail} should be taken with respect to
%$x$ as well.
%For both, continuous and differentiable cases,
%there is also an additional term
% \log[1-d_x^{d_x/2}\exp\{-R_m s_m\}/s_m^{d_x/2} ].

From the definition of models $\mathcal{M}_2$ and $\mathcal{M}_3$ and
from the comparison of the convergence rates
in (\ref{eq:bd_rate}) and (\ref{eq:bd_rate_linlogit}),
it is clear that using only linear indices in $x$ in the mixing
probabilities does not come without a cost.
The number of mixing components in model $\mathcal{M}_3$ that
approximates an infeasible model $\mathcal{M}_0$ is equal to $m N(m)+1$
while for model with
polynomial terms in logit, $\mathcal{M}_1$, this number is $m+1$
(Corollary \ref{crl:logit_polynomials}).
The proof of Corollary \ref{crl:linlogit_case_bounds} implies
that the number of hypercubes in the partition of $X$, $N(m)$,
increases exponentially with the dimensionality of $X$.
Thus, the number of parameters in model
$\mathcal{M}_3$ grows exponentially in
the dimension of $x$ (the exponential growth of the number of
parameters in $\mathcal{M}_1$ is discussed at the end of the previous section).
%Also,
%for both models smoothness of $F(A_j^m|x)$ in $x$ would somewhat
%mitigate this issue.
Overall, approximation results for $\mathcal{M}_1$ and $\mathcal{M}_3$
do not seem to suggest which model might perform better in practice;
however, they seem to identify a tradeoff between the number of
components in the mixture and the flexibility of models for the mixing
probabilities.

%s5 ###
\section{Flexible means and variances}
\label{sec:flex_mean}

In this section, I show that a finite mixture of normal regressions
models, in which
mixing probabilities do not depend on $x$, can be quite flexible.
However, the results also suggest that specifications in which
mixing probabilities are flexible functions of $x$ might perform better.

There is a large literature on finite mixture of regressions models.
In early work, mixtures of two normal regressions were considered [see,
e.g., \citet{QuandtRamsey78} and \citet{Kiefer78}].
\citet{JonesMcLachlan92} applied the EM algorithm for estimation
of finite mixtures of normal regressions.
Fitting of more general finite mixtures of generalized linear models
has been considered in \citet{Jansen93} and \citet
{WedelDeSarbo95} among others.
Many more references can be found in a comprehensive book on finite
mixture models by \citet{McLachlanPeel00}.

To the best of my knowledge, the literature on finite mixtures of
regressions does not contain any approximation results for conditional
densities.
The closest analogs of the results I obtain can be found in the
literature on finite mixtures of unconditional densities [see, e.g.,
\citet{ZeeviMeir1997} and references therein and \citet
{LiBarron99}].
Even for mixtures of unconditional densities
approximation results for the KL distance, which is useful for
establishing consistency of Bayesian or classical maximum likelihood
estimators, seem to be scarce.
Approximation results in the KL distance for convex combinations of
densities in \citet{ZeeviMeir1997} and \citet{LiBarron99}
seem to apply to mixtures of truncated normals and to target densities
that are compactly supported.
Some of these results are very strong.
For example, for target densities that are general mixtures of the
densities mixed in the model, approximation error bounds obtained by
\citet{LiBarron99}
are proportional to $m^{-1}$.
If there are no covariates $x$, then the infeasible model from Section
\ref{sec:Infeasible_model} is simply a finite mixture of multivariate normals.
For an elaboration on this idea in the context of joint and conditional
density estimation and for consistency results for a Bayesian estimator
based on this model see \citet{NoretsPelenis09}.
%(for an elaboration on this idea see \citet{NoretsPelenis09}).
The convergence rates obtained for this model in Section \ref
{sec:M0bounds} are slower than $m^{-1}$. However, the convergence rates
are not directly comparable as the target densities in \citet
{LiBarron99} are different from those considered here.

%results apply to more general settings, e.g., target densities and
%mixtures of densities with unbounded support.

%Actually, if there are no covariates $x$, then the infeasible model
%from Section \ref{sec:Infeasible_model} is simply a finite mixture of
%multivariate normals.
%For an elaboration on this idea in the context of joint and
%conditional density estimation and for consistency results for a
%Bayesian estimator based on this model see

Model $\mathcal{M}_4$ constructed in this section is very similar to
model $\mathcal{M}_0$ except for one important difference.
In $\mathcal{M}_4$,
fine equal probability partitions of $Y$ are used instead of
fine equal length partitions in $\mathcal{M}_0$.
As will be clear below, $\mathcal{M}_4$ defined in this way allows
mixing probabilities to be independent of $x$. However, it requires the
means of the mixed normals to be flexible functions of $x$.
In this section, I~assume that the response variable is univariate: $Y
\subset R$ or $d=1$ (all the results from previous sections were
obtained for arbitrary $d$).
If fine equal probability partitions
can be well defined for distributions of multivariate random variables
and if these partitions depend smoothly on covariates,
then it might be possible to extend
the results of this section to multivariate responses. I do not pursue
this conjecture here.

Define model $\mathcal{M}_4$ as follows:
\[
p(y|x ,\mathcal{M}_4) = \sum_{j=1}^m \alpha_j^m
\phi(y, \mu_j^m(x ),\sigma_j^m(x )).
\]
For a given $x$ let $A_j^m(x)$, $j=0,1,\ldots,m$, be a partition of $Y$
such that $\bigcup_{j=1}^m A_j^m(x)$ is a nondecreasing interval and
%
%e5.1 ###
%
\begin{eqnarray}
\label{eq:EPPconditions}
F(A_j^m(x)|x) &=& p_m,\qquad j>0,\nonumber\\[-8pt]\\[-8pt]
F(A_0^m(x)|x) &=& 1 - m p_m \quad\mbox{and}\quad m p_m
\rightarrow1,\nonumber
\end{eqnarray}
for some $p_m \in(0,m^{-1}]$ that does not depend on $x$.
Define an upper bound on the length of an element of the fine part of
the partition $h_m(x) \geq\break\max_{j>0} \lambda(A_j^m(x))$.
The candidate mixing probabilities are given by $\alpha
_j^m=F(A_j^m(x)|x)$ and
$\mu_j^m(x) \in A_j^m(x)$. The standard deviations
$\sigma_j^m(x )=\sigma_m(x)$ for $j>0$ and
$\sigma_0^m(x )=\sigma_0(x)$
are treated
as functions of $x$ which is not essential but it weakens the
restrictions on $F$ (Corollaries \ref{crl:flex_mean_bddsupport} and
\ref
{crl:poly_mean_bddsupport} and Examples \ref{ex:exponential_flex_mean}
and \ref{ex:uniform_flex_mean} below illustrate this point).
Note that $\mathcal{M}_4$ is an infeasible model;
in Corollary \ref{crl:poly_mean_bddsupport} below, I consider a
feasible model $\mathcal{M}_5$ in which $\mu_j^m(x)$ are approximated
by polynomials (see also Examples~\ref{ex:exponential_flex_mean}
and~\ref{ex:uniform_flex_mean}).

Suppose sequences $\delta_m(x)$, $\sigma_m(x)$, and $h_m(x)$ satisfy
%
%e5.2 ###
%
\begin{equation}
\label{eq:cond_delta_sigma_h_x}
\delta_m(x) \rightarrow0,\qquad \frac{\sigma_m(x)} {\delta_m(x)}
\rightarrow0, \qquad\frac{h_m(x)} {\sigma_m(x)} \rightarrow0.
\end{equation}
Next, let us introduce the following restrictions on $F$.
\begin{assumption} \label{assn:flexiblemu_case}
1.
Partitions $A_j^m(x)$ used in construction of $p(y|x ,\mathcal
{M}_4)$ satisfy
(\ref{eq:EPPconditions}), and (\ref{eq:cond_delta_sigma_h_x}) holds.
\smallskipamount=0pt
\begin{enumerate}[2.]
\item[2.]
\hypertarget{assnitem:flexiblemu_case_1}
$f(y|x)$ is continuous in $y$ a.s. $F$.%on $Y$ for all $x \in X$.
\item[3.]
\hypertarget{assnitem:flexiblemu_case_3}
For any $(y,x)$ there exists interval $C(r(x),y,x)$ with length
$r(x)>0$
and $y \in C(r(x),y,x)$
such that (i)
%
%e5.3 ###
%
\begin{equation}
\label{eq:IntBoundFinfFflexmu}
\int\log\frac{f(y|x)}{\inf_{z \in C(r(x),y,x)} f(z|x) } F(dy,dx) <
\infty
\end{equation}
and (ii) exists $M$ such that for any $m \geq M$,
if $y \in A_0^m(x)$, then
$C(r(x),y,x) \cap A_0^m(x)$
contains an interval $C_0(r(x),y,x)$ with an end at $y$ and length $r(x)/2$,
and
if $y \in Y \setminus A_0^m(x)$, then $C(r(x),y,x) \cap(Y \setminus A_0^m(x))$
contains
an interval $C_1(r(x),y,x)$ with an end at $y$ and length $r(x)/2$.
\item[4.] $h_m(x)$, $\sigma_m(x)$, and $r(x)$ satisfy
%
%e5.4 ###
%
\begin{equation}
\label{eq:cond_r_sigma_h_x}
\sup_x \frac{\sigma_m(x)} {r(x)} \rightarrow0,\qquad
\sup_x \frac{h_m(x) } {\sigma_m(x)} \rightarrow0.
\end{equation}
\item[5.]
$\sigma_0(x)$ and $r(x)$ satisfy
%
%e5.5 ###
%
\begin{equation}
\label{eq:cond_sigma0_x}
1> 1/4 \geq\phi(y, 0,\sigma_0(x)) r(x)/2,
\end{equation}
which holds, for example, when $\sigma_0(x) \geq2 (2 \pi)^{-1/2}
\cdot r(x)$.

\item[6.]
\hypertarget{assnitem:flexiblemu_case_2}
$| {\int\log[\phi(y, 0,\sigma_0(x)) r(x)/2 ] F(dy,dx)}| < \infty$.
%(More than the second moments of $y$ are finite).
\end{enumerate}
\end{assumption}
\begin{proposition} \label{prp:flexiblemu_case}
If Assumption \ref{assn:flexiblemu_case} holds then
$d_{\mathrm{KL}}(F,\mathcal{M}_4) \rightarrow0$ as \mbox{$m\rightarrow\infty$}.
\end{proposition}

The proposition is proved in the \hyperref[app]{Appendix}.
The assumptions of the proposition and their role in the proof are
similar to those discussed in detail in Section \ref
{sec:Infeasible_model} for~$\mathcal{M}_0$.
The assumptions are satisfied by a large class of densities as
illustrated by the following corollaries and examples.
Approximation error bounds for $\mathcal{M}_4$ are presented below in
Corollary \ref{crl:flex_mean_rate}.

\mbox{}
\begin{corollary} \label{crl:flex_mean_bddsupport}
Assume:
\begin{enumerate}
\item$f(y|x)$ is continuous in $y$ in the interior of the
support of $f(y|x)$ for all $x \in X$.
%on the support.$Y$ for all $x \in X$.
%
\item There exists $\overline{f}<\infty$, such that $ f(y|x) \leq
\overline{f}$ for all $(y,x)$.
\item The support of $f(\cdot|x)$ is given by a finite interval
$[a(x),b(x)]$, where $a(x)$ and $b(x)$ are square integrable.
Also, for some $\underline{f} \in(0,1)$, a positive integer $n$,
and\vspace*{1pt} $a(x) \leq a_1(x) \leq b_1(x) \leq b(x)$,%\\
$f(y|x) \geq\underline{f}$ on $[a_1(x) , b_1(x)]$, %\\
%
%f2 ###
%
\begin{figure}

\includegraphics{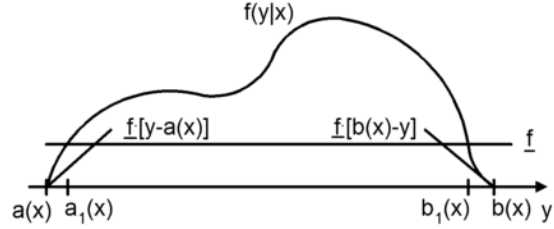}

\caption{Approximation of densities with bounded support by $\mathcal
{M}_4$.}%
\label{fig:fmubs}%
\end{figure}
$f(y|x) \geq\underline{f} \cdot[ y-a(x)]^n$ on $(a(x) ,
a_1(x))$, and %\\
$f(y|x) \geq\underline{f} \cdot[b(x)-y]^n$ on $(b_1(x) ,
b(x))$.
%Note that this implies $b(x)-a(x) \leq2 + \underline{f}^{-1}$.
Figure~\ref{fig:fmubs} provides an illustration for $n=1$.
\item
\hypertarget{crlitem:monotone}
There exists $r>0$ such that $f(\cdot|x)$ is nondecreasing on
$(a(x) , a_1(x)+r/2)$ and nonincreasing on
$(b_1(x) - r/2, b(x))$ for all $x \in X$.
\end{enumerate}
Then for
$\mathcal{M}_4$ constructed so that
$p_m=1/m$,
$A_0^m = \varnothing$,
%and $A_1^m(x),\ldots,A_m^m(x)$ is $F(\cdot|x)$ equal probability
%partition of $[a(x),b(x)]$,
$\mu_j^m(x) \in A_j^m(x)$ and
$\sigma_m(x)=p_m^{1/[4(n+1)]}$ and $\sigma_0(x) = 2 (2 \pi
)^{-1/2}\cdot
r$ are independent of $x$,
$d_{\mathrm{KL}}(F$,\break $\mathcal{M}_4) \rightarrow0$.
\end{corollary}
\begin{corollary} \label{crl:poly_mean_bddsupport}
Assume conditions from Corollary \ref{crl:flex_mean_bddsupport},
$F^{-1}(p|x)$ is continuous in $x$ for all $p \in[0,1]$, $X$ is
compact. Then there exists a sequence of polynomials $P_j^m(x )$ such that
$d_{\mathrm{KL}}(F,\mathcal{M}_5) \rightarrow0$ where
\[
p(y|x ,\mathcal{M}_5) = \sum_{j=1}^m p_m \phi(y, P_j^m(x ),p_m^{1/8}).
\]
\end{corollary}
\begin{pf}
Let $\mu_j^m(x)=F^{-1}((j-1/2)p_m|x)$.
Note that $\mu_j^m(x) \in A_j^m(x)=[F^{-1}((j-1)p_m|x),F^{-1}(j p_m|x)]$
and
\[
p_m/2 = \int_{\mu_j^m(x)}^{F^{-1}(j p_m|x)} f(y|x) \,dy \leq\bigl(F^{-1}(j
p_m|x) - \mu_j^m(x)\bigr) \overline{f}.
\]
Similarly,
$p_m/2 \leq(\mu_j^m(x) - F^{-1}((j-1) p_m|x)) \overline{f}$.
Thus, for $\varepsilon_m = p_m / (2 \overline{f})$,
$(\mu_j^m(x)-\varepsilon_m, \mu_j^m(x)+\varepsilon_m) \subset A_j^m(x)$.
By the Stone--Weierstrass theorem there exist finite order polynomials
in $x$, $P_j^m(x )$ such that
$|P_j^m(x ) - \mu_j^m(x)| < \varepsilon_m$. Therefore,
$P_j^m(x ) \in A_j^m(x)$, which was the only requirement on the means
of the mixture components in Corollary \ref{crl:flex_mean_bddsupport}.
\end{pf}
\begin{example}
\label{ex:exponential_flex_mean}
Exponential distribution, $f(y|x) = \gamma(x) \exp\{- \gamma(x)
y\}$,\break
$\gamma(x) \geq\underline{\gamma} > 0$, $\gamma(x)$ is continuous,
$\int\gamma \,dF < \infty$ and the second moment of $y$ is finite
($\int
\gamma^{-2} \,dF < \infty$).
The quantile function is given by $F^{-1}(p|x)= - \gamma(x)^{-1} \log(1-p)$.
Let the partition be such that $A_0^m = [F^{-1}(m p_m|x),\infty)$.
Since the exponential density is decreasing the largest interval in the
fine part of the partition is given by
$A_m^m = [ F^{-1}((m-1) p_m|x), F^{-1}( m p_m | x ) )$.
Therefore, $h_m(x)= h_m = \underline{\gamma}^{-1} \log(1+ p_m /
(1-p_m m))$.
Choosing $p_m = (m-m^{0.5})/m^2$ guarantees that
$h_m \rightarrow0$.
For $\sigma_m=
h_m^{1/4}$, and $\delta_m(x)=h_m^{1/8}$, and
$r(x)= 1$
conditions (\ref{eq:EPPconditions}), (\ref{eq:cond_delta_sigma_h_x})
and (\ref{eq:cond_r_sigma_h_x}) hold.

Next, let $C(1,y,x)=[y, y+1]$ if $y \in[0, 1/2]$,
$C(1,y,x)=[y-1/2, y+1/2]$ if $y \in[1/2, \infty)$.
Since
\[
\inf_{z \in C(1,y,x)} f(z|x) \geq\gamma(x) \exp\{-\gamma(x) (y+1)\},
\]
we have
\[
1 \leq f(y|x) \big/ \inf_{z \in C(1,y,x)} f(z|x) \leq\exp\{ \gamma(x) \}.
\]
Inequality (\ref{eq:IntBoundFinfFflexmu}) is satisfied since $\gamma
(x)$ is assumed to be integrable.
Finally, let $\sigma_0(x) = 2 (2 \pi)^{-1/2}$ so that
equation (\ref{eq:cond_sigma0_x}) in Assumption \ref
{assn:flexiblemu_case} holds.
Then,
\[
\biggl| \int\log[\phi(y, 0,\sigma_0(x)) r(x)/2 ] F(dy,dx)\biggr| =
\biggl| \int\biggl[-\log(4) - \frac{y^2 \pi}{4}\biggr] F(dy,dx) \biggr| < \infty
\]
since the second moment of $y$ is assumed to be finite.
Thus, condition \hyperlink{assnitem:flexiblemu_case_2}{6} of Assumption
\ref
{assn:flexiblemu_case} holds.

If $X$ is compact the same argument as in the proof of Corollary \ref
{crl:poly_mean_bddsupport} can be used to show that $\mu_j^m(x)$ can be
polynomial in $x$ [for fixed $m$ there exists $\varepsilon_m>0$ such that
$\lambda(A_j^m(x))>\varepsilon_m$ for all $x$ and $j$].

It is possible to give sufficient conditions for approximation results
when $\gamma(x)$ is not bounded away from zero, for example, let
$r(x)=\gamma(x)^{-1}$, $h_m(x) = \gamma(x)^{-1} \log(1+ p_m / (1-p_m
m))$, etc. However, then $\sigma_m$ and $\sigma_0$ would have to be
functions of $x$ [not necessarily flexible functions of $x$ but
functions that would have the same order as $\gamma(x)$]. Also,
$\gamma
(x)^{-1}$ is not continuous and the argument I use for justifying the
use of polynomial $\mu_j^m(x)$ breaks down in this case.
\end{example}
\begin{example}
\label{ex:uniform_flex_mean}
Uniform distribution, $f(y|x) = b(x)^{-1} 1_{[0,b(x)]}(y)$, $b(x)>0$ is
continuous,
$\int\log b \,dF < \infty$ and the second moment of $y$ is finite
($\int
b^{2} \,dF < \infty$).
This example demonstrates that the support of $f(y|x)$ does not have to
be (un)bounded uniformly in $x$ as long as normal variances are modeled
as flexible functions of $x$.

Let the\vspace*{-1pt} partition be such that $A_0^m = \varnothing$ and
$p_m=F(A_j^m|x)=m^{-1}$, $j>0$.
Note that $h_m(x)= b(x) / m$.
For\vspace*{1pt} $\sigma_m(x)=
b(x) p_m^{1/4}$, and $\delta_m(x)=b(x) p_m^{1/8}$, and
$r(x)= b(x)$
conditions (\ref{eq:EPPconditions}), (\ref{eq:cond_delta_sigma_h_x})
and (\ref{eq:cond_r_sigma_h_x}) hold.

Next, let $C(r(x),y,x)=[0,b(x)]$. Note that
$f(y|x) / \inf_{z \in C(r(x),y,x)} f(z|x) = 1$, and
inequality (\ref{eq:IntBoundFinfFflexmu}) is satisfied.
Finally, let $\sigma_0(x) = 2 (2 \pi)^{-1/2} b(x)$ so that
inequality (\ref{eq:cond_sigma0_x}) in Assumption \ref
{assn:flexiblemu_case} holds.
Then,
\[
\biggl| \int\log[\phi(y, 0,\sigma_0(x)) r(x)/2 ] F(dy,dx)\biggr| =
|{{-}\log(4) - \pi/ (3 \cdot4)}| < \infty
\]
and condition \hyperlink{assnitem:flexiblemu_case_2}{6} of Assumption
\ref
{assn:flexiblemu_case} holds.

If $X$ is compact and $b(x)$ is bounded away from zero then the same
argument, as in the proof of Corollary \ref{crl:poly_mean_bddsupport},
can be used to show that $\mu_j^m(x)$ can be polynomial in $x$ [for
fixed $m$ there exists $\varepsilon_m>0$ such that $\lambda
(A_j^m(x))>\varepsilon_m$ for all $x$ and $j$].
%All the arguments in the example are valid if expressions for
%$\sigma_0(x)$ and $\sigma_m(x)$ are defined with
%a function
%$\hat{b}(x)$ instead of $b(x)$ such that $\hat{b}(x)/b(x)$ is bounded
%above and away from zero.
\end{example}
\begin{corollary} \label{crl:flex_mean_rate}
Suppose conditions of Proposition
\ref{prp:flexiblemu_case} are satisfied for
$h_m(x)=h_m$, $\sigma_m(x)=\sigma_m$, $\delta_m(x)=\delta_m$ and
$r(x)=r$ that do not depend on $x$. Also, suppose
conditions from parts \textup{(i)} and \textup{(ii)} of Corollary \ref{crl:gen_case_bounds} hold.
Then for all sufficiently large $m$,
%
%e5.9 ###
%e5.8 ###
%e5.7 ###
%e5.6 ###
%
\begin{eqnarray}
\label{eq:bound1fm_intdsuplnFdz}
d_{\mathrm{KL}}(F,\mathcal{M}_4)
&\leq& \delta_m \cdot\frac{d^{1/2}}{2} \int\sup_{z \in C_{\delta
_m}(y)} \biggl\Vert\frac{d \log f(z|x) }{dz}\biggr\Vert F(dy,dx)
\\
\label{eq:bound2fm_lof1_epsm2}
&&{}
+ 2 \frac{3 h_m }{(2 \pi)^{1/2} \sigma_m}
+ 2 \exp\biggl\{-\frac{(\delta_m/\sigma_m)^2}{8} \biggr\}
\\
\label{eq:bound3fm_intdsuplnFdztail}
&&{} + \frac{r}{2} \int_{B_{\delta_m}(A_0^m(x))}
\sup_{z \in C_{r}(y)} \biggl\Vert\frac{d \log f(z|x) }{dz}\biggr\Vert F(dy,dx)
\\
\label{eq:bound4fm_inty2_cFtail_2}
&&{} + \int_{B_{\delta_m}(A_0^m(x))}
\biggl[ \frac{y^{\prime}y}{2 \sigma_0^2} -
\log\frac{(r/2)}{(2 \pi\sigma_0^2)^{1/2}} \biggr] F(dy,dx),
\end{eqnarray}
where
$B_{\delta_m}(A_0^m(x))=\{(y,x,)\dvtx C_{\delta_m}(y) \cap A_0^m(x) \neq
\varnothing\}
$%\]
and bounds in (\ref{eq:bound1fm_intdsuplnFdz})--(\ref
{eq:bound4fm_inty2_cFtail_2}) converge to zero as $m \rightarrow\infty$.
\end{corollary}
\begin{pf}
The proof is identical to the proof of Corollary \ref{crl:gen_case_bounds}.
\end{pf}

The bounds for $\mathcal{M}_4$, (\ref{eq:bound1fm_intdsuplnFdz})--(\ref{eq:bound4fm_inty2_cFtail_2}), are almost the same as
the bounds for $\mathcal{M}_0$, (\ref{eq:bound1_intdsuplnFdz})--(\ref{eq:bound4_inty2_cFtail_2}), obtained in Corollary \ref
{crl:gen_case_bounds}, except for a
difference between $B_{\delta_m}(A_0^m(x))$ in $\mathcal{M}_4$
and
$B_{\delta_m}(A_0^m)$ in $\mathcal{M}_0$.
For the same value of $h_m$,
the length of the complement of $A_0^m(x)$ in $\mathcal{M}_4$ is
bounded above by
$m h_m$ [$h_m=\max_{j>0} \lambda(A_j^m(x))$] which is the length of the
complement of $A_0^m$ in $\mathcal{M}_0$.
Thus the bounds obtained for $\mathcal{M}_4$ are likely to be larger
than the bounds obtained for $\mathcal{M}_0$.
Compact and interpretable conditions sufficient for deriving
an explicit approximation rate for $\mathcal{M}_4$ from (\ref
{eq:bound1fm_intdsuplnFdz})--(\ref{eq:bound4fm_inty2_cFtail_2}) seem to
be difficult to find. Instead, I show in the following example that not
only bounds for $\mathcal{M}_0$ can be smaller but also that
convergence for $\mathcal{M}_0$ can be slightly faster than for
$\mathcal{M}_4$.
\begin{example}
\label{ex:2side_exp_flex_mean_rates}
Laplace distribution, $f(y|x) = 0.5 \gamma(x) \exp\{- \gamma(x)
|y|\}$,\break
$\gamma(x) \geq\underline{\gamma} > 0$, $\gamma(x)$ is continuous,
$\int\gamma \,dF < \infty$ and the second moment of $y$ is finite
($\int
\gamma^{-2} \,dF < \infty$). Note that nondifferentiability of
$f(y|x)$ at zero does not affect any of the theoretical results above.

First\vspace*{-1pt} consider $\mathcal{M}_4$. Let $A_j^m(x)=[F^{-1}((1-p_m m)/2 +
(j-1)p_m|x), F^{-1}((1-p_m m)/2 + j p_m|x))$.
Note that $F^{-1}(p|x)=\log(2p)/\gamma(x)$ for $p<0.5$ and
$F^{-1}(p|x)=-\log(2(1-p))/\gamma(x)$ for $p \geq0.5$.
Then,
%
%e5.10 ###
%
\begin{eqnarray}
\label{eq:h_m_def}
h_m & \geq & F^{-1}\bigl((1-p_m m)/2 + p_m|x\bigr) - F^{-1}\bigl((1-p_m m)/2|x\bigr)
\nonumber\\[-8pt]\\[-8pt]
& = & \frac{1}{\gamma(x)} \log\biggl(1+ \frac{2 p_m}{1-p_m m} \biggr).
\nonumber
\end{eqnarray}
Since $h_m \rightarrow0$ and $m p_m \rightarrow1$ we can write
\[
p_m=\frac{1}{m+g(m)},
\]
where $g(m)$ satisfies $g(m)/m \rightarrow0$ and $g(m) \rightarrow
\infty$.
Note that
\[
B_{\delta_m}(A_0^m(x))
\subset
\biggl(-\infty, \frac{\log( 1-p_m m )(1-\varepsilon_0)}{\gamma
(x)} \biggr)
\cup\biggl(-\frac{\log( 1-p_m m )(1-\varepsilon_0)}{\gamma(x)},
\infty
\biggr)
\]
for any $\varepsilon_0 \in(0,1)$ and all sufficiently large $m$.
A direct calculation shows that
integrals in (\ref{eq:bound3fm_intdsuplnFdztail}) and
(\ref{eq:bound4fm_inty2_cFtail_2}) can be bounded by
\[
\mbox{constant} \cdot(1-p_m m)^{1-\varepsilon}
\leq\mbox{constant} \cdot\bigl(g(m)/m\bigr)^{1-\varepsilon}
\]
for any $\varepsilon\in(\varepsilon_0, 1)$ and all sufficiently large $m$.
From (\ref{eq:h_m_def}) and the mean value theorem,
\[
h_m \geq\mbox{constant} \cdot\gamma(x)^{-1} \cdot g(m)^{-1}.
\]
Since the approximation error bounds increase in $h_m$, we should
choose the smallest possible value for $h_m=\mbox{constant} \cdot
\underline{\gamma}^{-1} \cdot g(m)^{-1}$.
One can verify that
the smallest upper bound for
$\delta_m$, $h_m/\sigma_m$, $\exp\{-(\delta_m/\sigma_m)^2/8\}$ and
$(g(m)/m)^{1-\varepsilon}$ is inside the interval
$(m^{-1/3}, m^{-1/[3+\varepsilon_1]}]$
for any $\varepsilon_1>0$ and all sufficiently large $m$.
Thus,
\[
d_{\mathrm{KL}}(F,\mathcal{M}_4) \leq\mbox{constant} \cdot\biggl(\frac
{1}{m} \biggr)^{1/[3+\varepsilon_1]}.
\]

Next, consider $\mathcal{M}_0$.
Expressions (\ref{eq:bound3_intdsuplnFdztail})
and
(\ref{eq:bound4_inty2_cFtail_2}) are exponentially decreasing in $h_m m$.
Setting $h_m$ to a power of $m$,
one can show that
%Expressions in \eqref{eq:bound1_intdsuplnFdz} and
%Simplifying bounds \eqref{eq:bound1_intdsuplnFdz} -
%
\[
d_{\mathrm{KL}}(F,\mathcal{M}_0) \leq\mbox{constant} \cdot\biggl(\frac
{1}{m} \biggr)^{1/[2+\varepsilon_2]},
\]
for any $\varepsilon_2 > 0$ and all sufficiently large $m$.
These results suggest that $\mathcal{M}_0$ converges to the target
density faster than $\mathcal{M}_4$.
\end{example}

It might be unfair to compare approximation errors for $\mathcal{M}_0$
and $\mathcal{M}_4$. Although both models are ``infeasible'' and
include $m$ functions that need to be approximated by polynomials (or
splines), the error from approximation by the polynomials enters the
total approximation error in different ways. Nevertheless, the results
obtained in this section do seem to suggest that models in which mixing
probabilities depend on covariates might perform better in practice.

\section{\texorpdfstring{Comparison with Jiang and Tanner (\protect
\citeyear{JiangTanner99})}{Comparison with Jiang and Tanner
(1999)}}
\label{sec:comparison}

Jiang and Tanner (\citeyear{JiangTanner99}) is the only work on approximation of
conditional densities by ME that I am aware of. \citet{JiangTanner99}
develop approximation and estimation results for target densities
of the form
%
%e6.1 ###
%
\begin{equation}
\label{eq:exp_family}
\pi(y|x;h(\cdot)) = \exp\bigl(a(h(x))y + b(h(x)) + c(y)\bigr).
\end{equation}
Functions $a$, $b$ and $c$ are assumed to be known,
$a$ and $b$ are assumed to have nonzero derivatives and
$h(x)$ is assumed to have uniformly bounded continuous second order derivatives.
It seems that their results could still hold if
$a$, $b$ and $c$ are known only up to some parameters (see their Remark 4).
\citet{JiangTanner99} show that $\pi(y|x;h(\cdot))$ can be
approximated in the KL distance by ME of the form
%
%e6.2 ###
%
\begin{equation}
\label{eq:JTmodel}
\sum_{j=1}^m \alpha_j^m(x) \pi(y|x;h_j(\cdot)),
\end{equation}
where $\pi(\cdot| \cdot;\cdot)$ is defined in (\ref{eq:exp_family}),
$h_j(x)$ is a linear function of $x$ and the mixing probabilities
$\alpha_j^m(x)$ can be modeled by logit (more general specifications
for mixing weights are also allowed).
The idea of their argument is to divide $X$ into a fine partition
$B_j^m$, approximate $1_{B_j^m}(x)$ by $\alpha_j^m(x)$ and approximate
$h(x)$ by linear function $h_j(x)$ on $B_j^m$.
\citet{JiangTanner99} prove that for their target class of densities
a bound on the approximation error is proportional to $m^{-4/d_x}$.

There are several important differences between the present work and
\citet{JiangTanner99}.
First, I consider multivariate responses, $y$, while \citet
{JiangTanner99} consider univariate responses.
Most importantly, I do not assume that functional form of $f(y|x)$ is
known, for example, known $\pi$, $a$, $b$ and $c$.
%Instead, I use
%weak regularity conditions such as restrictions on the support,
%smoothness and tail behavior to describe densities that can be
%approximated.
The components of the model I employ, for example, normal densities and
logit mixing probabilities, are generally not related to the true density.
As Examples \ref{ex:inf_student} and \ref{ex:inf_cont_bddsupport} and
Corollary \ref{crl:flex_mean_bddsupport} illustrate,
many densities that are not from (\ref{eq:exp_family}) are shown to be
approximable by ME models.
Examples \ref{ex:exponential} and \ref{ex:exponential_flex_mean} also
show that some of the densities from class (\ref{eq:exp_family})
satisfy sufficient conditions for approximation results I obtain.
However, there might exist densities from (\ref{eq:exp_family}) that
violate these sufficient conditions. This would not be surprising since
the ``correct'' functional forms are mixed in (\ref{eq:JTmodel}).
For the same reason it is not surprising that the approximation rate
obtained by \citet{JiangTanner99}, $m^{-4/d_x}$, differs from the
ones obtained here, for example, $m^{-1/[d_x+2+1/(q-2)+\varepsilon]}$ for
model $\mathcal{M}_3$ in Corollary \ref{crl:linlogit_case_bounds}.
% for specification $\mathcal{M}_3$ and $m^{-1/[2+1/(q-2)+\varepsilon]}$
%for specification $\mathcal{M}_5$. It is worth emphasizing that
%the approximation rate for $\mathcal{M}_5$ does not depend on the
%dimensions of $y$ and $x$.

%An advantage of the results in \citet{JiangTanner99} is that they
%obtain
%approximation rate as a simple function of $m$: a constant times .
%Bounds on the approximation error I obtain in this paper do not have
%an explicit expression that is proportional to a function of $m$ only.
%However, they do provide insight into how
%smoothness and tail behavior of the target density affect
%approximation error as $m$ increases. They also show that the
%dimensions of $x$ and $y$ enter the approximation error bounds
%exponentially.

Finally, responses in \citet{JiangTanner99} class (\ref
{eq:exp_family}) can be discrete, for example, Poisson. To accommodate
discrete responses in the framework of the present paper one could
map the discrete values of response $y$ into a partition of $R$
and introduce a corresponding latent variable $y^* \sim
p(y^*|x,\mathcal{M})$.
For example, for binary $y \in\{0,1\}$ let $y^* \in(-\infty,0)$ if
$y=0$ and $y^* \in[0,\infty)$ if $y=1$.
Any discrete distribution can be represented by a continuously
distributed latent variable in this fashion. This continuous
distribution can be flexibly modeled by $p(y^*|x,\mathcal{M})$.
Models with latent variables are easy to estimate in the Bayesian
framework using MCMC methods [see, e.g., \citet{TannerWong87} and
\citet{AlbertChibbinpoly93}].

%s7 ###
\section{Discussion}
\label{sec:conclusion}

This paper shows that large classes of conditional densities can be
approximated in the Kullback--Leibler distance by different
specifications of finite smooth mixtures of normal densities or
regressions. The theory can be generalized to smooth mixtures of
location scale densities. These results have interesting implications
for applied researchers.

First of all, smooth mixtures of densities or experts can be used as
flexible models for estimation of multivariate conditional densities.
It seems this issue has not been explored in the literature and it
would be interesting to see how specifications studied in the paper
work in these settings.

Second, smooth mixtures of simple components, for example, models in
which mixing probabilities are modeled by multinomial logit linear in
covariates and the means and variances do not depend on covariates, can
be quite flexible.
A~simulation study in \citet{VillaniKohnGiordani07} suggests
though that
models with more complex components perform better in practice.
This issue should be further explored in simulation studies.

Third, results in Section \ref{sec:lin_logit} suggest that making
mixing probabilities more flexible, for example, by using polynomials
in logit, might reduce the number of necessary mixture components.
However, these models are more difficult to estimate.

Fourth, models in which mixing probabilities do not depend on
covariates can be very flexible at least for univariate response
variables. However, they seem to require a lot of mixture components and
very flexible models for the means of the mixed normals. Also,
approximation error bounds and convergences rates (Example \ref
{ex:2side_exp_flex_mean_rates}) obtained in Section \ref{sec:flex_mean}
suggest that models with flexible mixing probabilities might perform
better in practice than models with flexible means of the mixed normals
and constant mixing probabilities.
Nevertheless, it would be interesting to see how these specifications
perform in actual applications and simulation studies.

On the basis of a simulation study,
\citet{VillaniKohnGiordani07} generally recommend using
heteroscedastic experts (mixture components with variances that depend
on covariates).
The theory obtained here suggests that heteroscedastic experts might be
necessary when differences in quantiles of $f(\cdot|x)$ are not
uniformly bounded in $x$ and, especially,
when the support bounds of $f(\cdot|x)$ are increasing without a bound
in $x$ (see Examples \ref{ex:infeasible_unifrm} and
\ref{ex:uniform_flex_mean}).
This suggestion is likely to remain useful
when the differences in quantiles and/or support of $f(\cdot|x)$,
although bounded, still change considerably with covariates.

%Finally, for models with the mixing probabilities modeled by scaled
%and normalized Gaussian densities, the approximation rates are shown
%to be independent of dimensions of $y$ and $x$. For all other
%specifications considered in the paper obtained approximation rates or
%the number of parameters are affected exponentially by the dimension
%of $x$.

Practical implications of the theoretical results obtained in the paper
and summarized in this section are deduced under the assumption of no
estimation and parameter uncertainty.
Exploring the behavior of the estimation error in addition to the
approximation error would result in a more complete understanding of
the ME models. This issue is left for future work.

Overall, the paper provides a number of encouraging approximation
results for (smooth) mixtures of densities or experts which might
stimulate more theoretical and applied work in this area of research.

%s8 ###
%
\begin{appendix}
\section*{Appendix}\label{app}

\vspace*{-14pt}

\begin{pf*}{Proof of Proposition \protect\ref{prp:gen_case}}
Since $d_{\mathrm{KL}}$ is always nonnegative,
\[
0 \leq\int\log\frac{f(y|x)}{p(y|x ,\mathcal{M}_0)} F(dy,dx) \leq
\int\log\max\biggl\{1, \frac{f(y|x)}{p(y|x ,\mathcal{M}_0)} \biggr\} F(dy,dx).
\]
Thus, it suffices to show that the last integral in the inequality
above converges to zero as $m$ increases. The dominated convergence
theorem (DCT) is used for that.
First, I establish conditions for point-wise convergence of the
integrand to zero a.s. $F$. Then, I present conditions for existence of
an integrable upper bound on the integrand required by the DCT.

For fixed $(y,x)$,
%
%e8.1 ###
%
\begin{eqnarray}\label{eq:model_lb}
p(y|x ,\mathcal{M}_0) &=&
\sum_{j=1}^m F(A_j^m|x) \phi(y, \mu_j^m,\sigma_m)
+F(A_0^m|x) \phi(y, 0,\sigma_0) \nonumber\\[-8pt]\\[-8pt]
&\geq&
\inf_{z \in C_{\delta_m}(y)} f(z|x)
\sum_{j\dvtx A_j^m \subset C_{\delta_m}(y)} \lambda(A_j^m) \phi(y, \mu
_j^m,\sigma_m),\nonumber
\end{eqnarray}
where $\lambda$ is the Lebesgue measure.

In Lemmas \ref{lm:boundRSbyInt_ESP} and \ref{lm:boundIntBy1}, I derive
the following bounds
for the Riemann sum in (\ref{eq:model_lb})
(the Riemann sum is not far from the corresponding normal integral, and
the integral is not far from 1):
%
%e8.2 ###
%
\begin{eqnarray}
\label{eq:RiemannSumBd}
&& \sum_{j\dvtx A_j^m \subset C_{\delta_m}(y)} \lambda(A_j^m) \phi(y,
\mu
_j^m,\sigma_m) \nonumber\\
&&\qquad \geq1 - \frac{3 d^{3/2} \delta_m^{d-1} h_m }{(2 \pi)^{d/2} \sigma
_m^d} -
\frac{8 (\sigma_m / \delta_m) }{(2 \pi)^{1/2} }
\exp\biggl\{-\frac{(\delta_m/\sigma_m)^2}{8} \biggr\}
\\
&&\qquad \geq1 - \frac{3 d^{3/2} \delta_m^{d-1} h_m }{(2 \pi)^{d/2} \sigma
_m^d} -
\exp\biggl\{-\frac{(\delta_m/\sigma_m)^2}{8} \biggr\},
\nonumber
\end{eqnarray}
where the last inequality holds for all sufficiently large $m$ ($\delta
_m/\sigma_m \rightarrow\infty$).
Given $\varepsilon> 0$ there exists $M_1$ such that for $m \geq M_1$,
expressions in
(\ref{eq:RiemannSumBd}) are bounded below by $(1 - \varepsilon)$.

If $f(y|x)$ is continuous in $y$ at
$(y,x)$ and $f(y|x)>0$
there exists $M_2$ such that for $m \geq M_2$,
$[f(y|x) / \inf_{z \in C_{\delta_m}(y)} f(z|x)] \leq(1 + \varepsilon)$
since $\delta_m \rightarrow0$.
For any $m \geq\max\{M_0, M_1, M_2\}$,
\begin{eqnarray*}
1 &\leq& \max
\biggl\{
1,
\frac{f(y|x)}{p(y|x ,\mathcal{M}_0)} \biggr\}
\\
&\leq&
\max\biggl\{1,\frac{f(y|x)}{\inf_{z \in C_{\delta_m}(y)} f(z|x) (1 -
\varepsilon
) } \biggr\}\leq\frac{1 + \varepsilon}{1 - \varepsilon}.
\end{eqnarray*}
Thus, $\log\max\{1,f(y|x)/p(y|x ,\mathcal{M}_0)\} \rightarrow0$
a.s. $F$ as long as $f(y|x)$ is continuous in $y$ a.s. $F$
[$f(y|x)$ is always positive a.s. $F$].

Parts \hyperlink{assnitem:general_case_2}{2} and \hyperlink
{assnitem:general_case_3}{3}
of Assumption \ref{assn:general_case} are used for establishing an
integrable upper bound for the DCT
%
%e8.3 ###
%
\begin{eqnarray}\label{eq:RiemannSumInBound}
p(y|x ,\mathcal{M}_0) &=&
\sum_{j=1}^m F(A_j^m|x) \phi(y, \mu_j^m,\sigma_m)
+F(A_0^m|x) \phi(y, 0,\sigma_0) \nonumber\\
&\geq&
[1-1_{A_0^m}(y)] \nonumber\\[-8pt]\\[-8pt]
&&\hspace*{0pt}{}\times\inf_{z \in C_1(r,y,x)} f(z|x) \cdot
\sum_{j\dvtx A_j^m \subset C_1(r,y,x)} \lambda(A_j^m) \phi(y, \mu
_j^m,\sigma_m)
\nonumber\\
&&{} + 1_{A_0^m}(y) \cdot
\inf_{z \in C_0(r,y,x)} f(z|x) \cdot
\lambda(C_0(r,y,x)) \phi(y, 0,\sigma_0). \nonumber
\end{eqnarray}
Lemmas \ref{lm:boundRSbyInt_ESP} and \ref{lm:boundIntBy1} imply that the
Riemann sum in (\ref{eq:RiemannSumInBound}) is bounded
below by $2^{-d} - 2^{-(d+1)}=2^{-(d+1)}$ for any $m$ larger then some $M_4$.
Inequalities (\ref{eq:RiemannSumInBound}) and (\ref{eq:cond_sigma0}) imply
%
%e8.4 ###
%
\begin{eqnarray}\label{eq:inqlty1}\qquad
&&\log\max\biggl\{1,\frac{f(y|x)}{p(y|x ,\mathcal{M}_0)}\biggr\}
\nonumber\\
&&\qquad\leq
\log\max\biggl\{1,\frac{f(y|x)}{\inf_{z \in C(r,y,x)} f(z|x) \cdot\phi(y,
0,\sigma_0) \cdot(r/2)^d}\biggr\} \nonumber\\[-8pt]\\[-8pt]
&&\qquad=
\log\frac{1}{\phi(y, 0,\sigma_0) (r/2)^d}
\max\biggl\{\phi(y, 0,\sigma_0) (r/2)^d,\frac{f(y|x)}{\inf_{z \in C(r,y,x)}
f(z|x)}\biggr\} \nonumber\\
&&\qquad\leq
-\log\bigl(\phi(y, 0,\sigma_0) (r/2)^d\bigr) +
\log\frac{f(y|x)}{\inf_{z \in C(r,y,x)} f(z|x)},\nonumber
\end{eqnarray}
where inequality (\ref{eq:inqlty1}) follows by the first inequality in
(\ref{eq:cond_sigma0}).
The first expression in (\ref{eq:inqlty1})
is integrable by
Assumption \ref{assn:general_case}, part \hyperlink{assnitem:general_case_2}{2}.
The second expression in (\ref{eq:inqlty1}) is integrable by
Assumption \ref{assn:general_case}, part \hyperlink{assnitem:general_case_3}{3}.
%This finishes the proof of the following proposition.
Thus the proposition is proved.
\end{pf*}
\begin{pf*}{Proof of Corollary \protect\ref{crl:gen_case_bounds}}
The proof of the first part of the proposition is a simple implication
of the argument in the proof of Proposition \ref{prp:gen_case}.
Note that
%
%e8.5 ###
%
\begin{eqnarray}
\label{eq:d_kl_2_ints}
d_{\mathrm{KL}}(F,\mathcal{M}_0) &=& \int_{Y\times X \setminus B_{\delta_m}(A_0^m)}
\log\frac{f(y|x)}{p(y|x,\mathcal{M}_0)} F(dy,dx) \nonumber\\[-8pt]\\[-8pt]
&&{} +
\int_{B_{\delta_m}(A_0^m)}
\log\frac{f(y|x)}{p(y|x,\mathcal{M}_0)} F(dy,dx).
\nonumber
\end{eqnarray}
For $(y,x) \in Y\times X \setminus B_{\delta_m}(A_0^m)$,
inequalities (\ref{eq:model_lb}) and (\ref{eq:RiemannSumBd})
apply.
Thus, the first integral in (\ref{eq:d_kl_2_ints}) is bounded by the
sum of
(\ref{eq:bound1_intFinfF})
and (\ref{eq:bound2_lof1_epsm}), where the bound in
(\ref{eq:bound2_lof1_epsm}) is obtained
by the mean value theorem for $-\log(1-x)$ and a small positive $x$,
%
%e8.6 ###
%
\begin{eqnarray}
\label{eq:bd_rate_log}
&& - \log\biggl(
1 - \frac{3 d^{3/2} \delta_m^{d-1} h_m }{(2 \pi)^{d/2} \sigma_m^d} -
\exp\biggl\{-\frac{(\delta_m/\sigma_m)^2}{8} \biggr\}
\biggr) \nonumber\\[-8pt]\\[-8pt]
&&\qquad \leq2 \biggl(
\frac{3 d^{3/2} \delta_m^{d-1} h_m }{(2 \pi)^{d/2} \sigma_m^d} +
\exp\biggl\{-\frac{(\delta_m/\sigma_m)^2}{8} \biggr\}
\biggr).
\nonumber
\end{eqnarray}
By inequality (\ref{eq:RiemannSumInBound}),
the second integral in (\ref{eq:d_kl_2_ints}) is bounded by the sum of
(\ref{eq:bound3_intFinfFtail})
and (\ref{eq:bound4_inty2_cFtail}).

Expression (\ref{eq:bound1_intFinfF}) converges to zero by the DCT. The
point-wise convergence follows by the assumed continuity and positivity
of $f(y|x)$. An integrable upper bound is given by (\ref{eq:IntBoundFinfF}).
Expression (\ref{eq:bound1_intFinfF}) converges to zero by
(\ref{eq:cond_delta_sigma_h}).
Expressions (\ref{eq:bound3_intFinfFtail}) and (\ref
{eq:bound4_inty2_cFtail}) converge to zero because
$Y\times X \setminus B_{\delta_m}(A_0^m) \nearrow Y\times X$ and the
integrands are integrable by (\ref{eq:IntBoundFinfF}) and by the
assumed finiteness of the second moment of $y$.
Thus, the first part of the proposition is proved.

The second part of the proposition [bounds for differentiable $f(y|x)$]
follows from the first part since
\[
\biggl| \log\frac{f(y|x)}{\inf_{z \in C_{r}(y)} f(z|x)} \biggr| \leq
\sup_{z \in C_{r}(y)} \biggl\Vert\frac{d \log f(z|x) }{dz}\biggr\Vert\frac
{d^{1/2} r}{2},
\]
which is implied by
the multivariate mean value theorem: for any $(z_1,z_2)$
\[
|{\log f(z_1|x) - \log f(z_2|x) }| \leq\bigl\Vert f^{\prime}\bigl(c z_1 + (1-c) z_2\bigr)
\bigr\Vert\Vert z_1-z_2\Vert
\]
for some $c \in[0,1]$.
Convergence of the bounds to zero is obtained in the same way as in the
first part of the proposition.

To obtain the third part let us suppose that the fine part of the
partition $\{A_j^m, 1 \leq j \leq m \} $ is
centered at 0.
If $(y,x) \in B_{\delta_m}(A_0^m)$, then $|y_i| \geq h_m m^{1/d} / 2 -
\delta_m > h_m m^{1/d}/3$ for $i \in\{1,\ldots,d\}$ and all
sufficiently large $m$ and
%
%e8.7 ###
%
\begin{eqnarray} \label{eq:bd_rate_y2}\qquad
&& \int_{B_{\delta_m}(A_0^m)}
y_i^2 F(dy,dx) \nonumber\\
&&\qquad\leq
\int_{\{(y,x)\dvtx|y_i| > h_m m^{1/d} / 3, \forall i \}}
y_i^2 F(dy,dx)\nonumber\\[-18pt] \\
&&\qquad \leq
(h_m m^{1/d} / 3)^{-(q-2)}\nonumber\\
&&\qquad\quad{}\times \int_{ \{(y,x)\dvtx|y_i| > h_m m^{1/d} / 3,
\forall i \}}
(h_m m^{1/d} / 3)^{q-2} y_i^2 F(dy,dx) \nonumber\\
&&\qquad
\leq
(h_m m^{1/d} / 3)^{-(q-2)} \int_{ Y \times X}
y_i^q F(dy,dx).\nonumber
\end{eqnarray}
Similarly,
%
%e8.8 ###
%
\begin{eqnarray}\label{eq:bd_rate_dlogfdz}\qquad
&& \int_{B_{\delta_m}(A_0^m)}
\sup_{z \in C_{r}(y)} \biggl\Vert\frac{d \log f(z|x) }{dz}\biggr\Vert F(dy,dx)
\nonumber\\
&&\qquad
\leq
\int_{\{(y,x)\dvtx|y_i| > h_m m^{1/d} / 3, \forall i \}}
\sup_{z \in C_{r}(y)} \biggl\Vert\frac{d \log f(z|x) }{dz}\biggr\Vert F(dy,dx)
\nonumber
\\
&&\qquad \leq
\biggl(\int_{ \{(y,x)\dvtx|y_i| > h_m m^{1/d} / 3, \forall i \}}
(h_m m^{1/d} / 3)^{q-2}\nonumber\\[-8pt]\\[-8pt]
&&\hspace*{137pt}{}\times \sup_{z \in C_{r}(y)} \biggl\Vert\frac{d \log f(z|x)
}{dz}\biggr\Vert F(dy,dx)\biggr)\nonumber\\
&&\qquad\quad\hspace*{0pt}{}\times\bigl((h_m m^{1/d} / 3)^{q-2} \bigr)^{-1} \nonumber\\
&&\qquad
\leq
(h_m m^{1/d} / 3)^{-(q-2)} \int_{ Y \times X}
y_{i_1}^{q-2} \sup_{z \in C_{r}(y)} \biggl\Vert\frac{d \log f(z|x)
}{dz}\biggr\Vert F(dy,dx).\nonumber
\end{eqnarray}

%Let $\delta_m = m^{-\alpha_{\delta}}$, $\sigma_m = m^{-\alpha_{
%$h_m = m^{-\alpha_{h}}$.

Since integrals in (\ref{eq:bd_rate_y2}) and (\ref{eq:bd_rate_dlogfdz})
are finite by assumption, (\ref{eq:bound3_intdsuplnFdztail}) and (\ref
{eq:bound4_inty2_cFtail_2}) can be bounded above by an expression
proportional to $(h_m m^{1/d})^{-(q-2)}$.
Thus, the sum of (\ref{eq:bound1_intdsuplnFdz})--(\ref
{eq:bound4_inty2_cFtail_2}) is bounded by
%
%e8.9 ###
%
\begin{eqnarray}
\label{eq:bd_rate_alphas}
&&c_1 \cdot\delta_m
+c_2 \cdot\exp\{- (\delta_m/\sigma_m)^2/8\}
+c_3 \cdot\delta_m^{d-1} h_m / \sigma_m^d\nonumber\\[-8pt]\\[-8pt]
&&\qquad{}+c_4 \cdot1/(h_m m^{1/d})^{q-2},\nonumber
\end{eqnarray}
where constants $c_1$, $c_2$, $c_3$ and $c_4$ do not depend on $m$.
Let $b_m$ be the smallest number satisfying
$b_m \geq\delta_m$,
$b_m \geq\delta_m^{d-1} h_m / \sigma_m^d$, $b_m \geq1/(h_m
m^{1/d})^{q-2}$ and $b_m \geq\exp\{- (\delta_m/\sigma_m)^2/8\}$.
The first three of these inequalities imply
\[
b_m \geq\{[(\delta_m/\sigma_m)^d]/m^{1/d}\}^{1/[2+1/(q-2)]}.
\]
It implies that for all sequences $\delta_m$, $\sigma_m$ and $h_m$
allowed by the corollary,
\[
b_m > \biggl( \frac{1}{m} \biggr)^{1/(d\cdot[2+1/(q-2)])}.
\]
One can verify that
%
%e8.10 ###
%
\begin{equation}\qquad
\label{eq:2bds_rate}
b_m \leq\biggl( \frac{(4 \log m / d )^{d/2}}{m^{1/d}} \biggr)^{1/[2+1/(q-2)]}
\leq\biggl( \frac{1}{m} \biggr)^{1/(d \cdot[2+1/(q-2)+\varepsilon])},
\end{equation}
when $\delta_m$ equal to the first bound in (\ref{eq:2bds_rate}),
$(\delta_m/\sigma_m)^2 = 4 \log m /d$ and $h_m=\delta_m^2/(\delta
_m/\sigma_m)^d$.
\end{pf*}
\begin{pf*}{Proof of Proposition \protect\ref{prp:gen_linear_logit}}
Define $I_1^m(x,s_m)=\{i\dvtx\Vert x_i^m-x\Vert^2 < s_m \}$ and
$I_2^m(x,s_m)=\{
i\dvtx\Vert x_i^m-x\Vert^2 > 2 s_m \}$.
For $i \in I_1^m(x,s_m)$,
%
%e8.11 ###
%
\begin{equation}
\label{eq:I1Rmsm}
[- R_m (x_i^{m \prime} x_i^m - 2 x_i^{m \prime} x )] > [-R_m (s_m -
x^{ \prime} x)]
\end{equation}
and for $i \in I_2^m(x,s_m)$,
%
%e8.12 ###
%
\begin{equation}
\label{eq:I2Rmsm}
[- R_m (x_i^{m \prime} x_i^m - 2 x_i^{m \prime} x )] < [-R_m (2 s_m -
x^{ \prime} x)].
\end{equation}
Note that
%
%e8.13 ###
%
\begin{eqnarray}
\label{eq:BoundSumI1}\qquad
&&\frac{ \sum_{i \in I_1^m(x,s_m)}
\exp\{ - R_m (x_i^{m \prime} x_i^m - 2 x_i^{m \prime} x ) \}}
{\sum_{l} \exp\{ - R_m (x_l^{m \prime} x_l^m - 2 x_l^{m \prime} x
)\} }
\nonumber\\
&&\qquad \geq
1- \frac{ \sum_{i \in I_2^m(x,s_m)}
\exp\{ - R_m (x_i^{m \prime} x_i^m - 2 x_i^{m \prime} x ) \}}
{\sum_{i \in I_1^m(x,s_m)} \exp\{ - R_m (x_i^{m \prime} x_i^m - 2
x_i^{m \prime} x )\} }
\\
&&\qquad
\geq
1 -
\frac{\mbox{card}(I_2^m(x,s_m))}{\mbox{card}(I_1^m(x,s_m))}
\exp\{-R_m s_m\}
\geq1 - d_x^{d_x/2} \frac{\exp\{-R_m s_m\}}{s_m^{d_x/2}},
\nonumber
\end{eqnarray}
where the second inequality follows from (\ref{eq:I1Rmsm}) and (\ref
{eq:I2Rmsm}). The last inequality
follows from the following
bounds on the number of elements in $I_1^m(x,s_m)$ and $I_2^m(x,s_m)$:
$\mbox{card}(I_1^m(x,s_m)) \geq1$
[$s_m$ is chosen in (\ref{eq:CondSmRm}) so that any ball in $X$ with
radius $s_m^{1/2}$ has to contain at least one
$x_i^m$] and
\[
\mbox{card}(I_2^m(x,s_m)) \leq N(m) = d_x^{d_x/2} s_m^{-d_x/2}.
\]
For $i \in I_1^m(x,s_m)$ and $A_j^m \subset C_{\delta_m}(y)$,
%
%e8.14 ###
%
\begin{equation}
\label{eq:FA_jBoindLogitLin}
F(A_j^m|x_i^m) \geq\lambda(A_j^m) \inf_{z \in C_{\delta_m}(y),
\Vert t-x\Vert^2 \leq s_m} f(z|t).
\end{equation}
Inequalities (\ref{eq:BoundSumI1}), (\ref{eq:FA_jBoindLogitLin}) and
(\ref{eq:RiemannSumBd}) imply that $p(y|x ,\mathcal{M}_3)$ exceeds
\begin{eqnarray*}
&& \sum_{j\dvtx A_j^m \subset C_{\delta_m}(y)} \sum_{i \in I_1^m(x,s_m)}
F(A_j^m|x_i^m)
\frac{
\exp\{ - R_m (x_i^{m \prime} x_i^m - 2 x_i^{m \prime} x ) \}}
{\sum_{l} \exp\{ - R_m (x_l^{m \prime} x_l^m - 2 x_l^{m \prime} x
)\} }
\phi(y, \mu_j^m,\sigma_m) \\
&&\qquad \geq
\inf_{z \in C_{\delta_m}(y), \Vert t-x\Vert^2 \leq s_m} f(z|t) \cdot
\biggl[
1 - \frac{3 d^{3/2} \delta_m^{d-1} h_m }{(2 \pi)^{d/2} \sigma_m^d}\\
&&\qquad\quad\hspace*{120.2pt}{} -
\frac{8 d \sigma_m }{(2 \pi)^{1/2} \delta_m}
\exp\biggl\{-\frac{(\delta_m/\sigma_m)^2}{8} \biggr\}
\biggr]
\\
&&\qquad\quad\hspace*{78.5pt}{} \times
\biggl[1-d_x^{d_x/2}\frac{\exp\{-R_m s_m\}}{s_m^{d_x/2}} \biggr].
\end{eqnarray*}
The expression on the last line of this inequality converges to $1$ by
(\ref{eq:CondSmRm}).
The rest of the proof is exactly the same as
%follows by the same argument as in
the proof of Proposition \ref{prp:gen_case}.
%As long as
%conditions in \eqref{eq:CondSmRm} and Proposition \ref{prp:gen_case}
%hold,
%for any $(y,x)$ such that
%$f(y|x)$ is continuous in $(y,x)$ and $f(y|x)>0$,
%which implies convergence a.s. $F$ by part
%Assumption \ref{assn:logitlin_case}.
%Assumption \ref{assn:logitlin_case} also implies existence of an
%integrable upper bound
%in exactly the same way Assumption \ref{assn:general_case} does for $
%The proposition claim follows by the DCT.
\end{pf*}
\begin{pf*}{Proof of Corollary \protect\ref{crl:linlogit_case_bounds}}
The proof of part (i) is identical to the proof of Corollary \ref
{crl:gen_case_bounds} part (ii).

The proof of part (ii) is also similar to the proof of Corollary \ref
{crl:gen_case_bounds} part (iii). Just set $s_m^{1/2} = \delta_m$ and
note that (\ref{eq:bound5_expRmsm}) can be made arbitrarily smaller
than the other parts of the bound by an appropriate choice of $R_m$.
Thus, the bound is the same as in (\ref{eq:bd_rate}), we just\vspace*{1pt} need to
express $m$ in terms of the number of mixture components in $\mathcal
{M}_3$, $m N(m)$.
From the definition of $N(m)$ and $s_m$,
$N(m)=\lambda(B_i^m)^{-1}=d_x^{d_x/2} s_m^{-d_x/2}$.
Since we set $s_m^{1/2} = \delta_m$ and $\delta_m=m^{-1/(d\cdot
[2+1/(q-2)])}$ in the proof of Corollary \ref{crl:gen_case_bounds},
\[
m N(m)= d_x^{d_x/2} m^{1+d_x / (d\cdot[2+1/(q-2)])}.
\]
From this equation, one can express $m$ as a function of
$m N(m)$ and plug it in (\ref{eq:bd_rate}) to obtain
(\ref{eq:bd_rate_linlogit}).
\end{pf*}
\begin{pf*}{Proof of Proposition \protect\ref{prp:flexiblemu_case}}
First, consider point-wise convergence a.s. $F$.
For fixed $(y,x)$ and an interval $C_{\delta_m(x)}(y)$ with center $y$
and length $\delta_m(x) > 0$,
%
%e8.15 ###
%
\begin{eqnarray} \label{eq:M4_lb}\hspace*{32pt}
p(y|x ,\mathcal{M}_4) &=&
\sum_{j=1}^m F(A_j^m(x)|x) \phi(y, \mu_j^m(x),\sigma_m(x))
\nonumber\\
&&{}+F(A_0^m(x)|x) \phi(y, 0,\sigma_0(x)) \nonumber\\
&\geq&
\inf_{z \in C_{\delta_m(x)}(y)} f(z|x)
\sum_{j=1}^m \lambda\bigl(A_j^m(x) \cap
C_{\delta_m(x)}(y)\bigr)\nonumber\\[-8pt]\\[-8pt]
&&\hspace*{90.1pt}{}\times \phi(y, \mu
_j^m(x),\sigma_m(x)) \nonumber\\
&\geq&
\inf_{z \in C_{\delta_m(x)}(y)} f(z|x)
\biggl(1 -
\frac{ 6 h_m(x) } {(2 \pi)^{1/2} \sigma_m(x) }\nonumber\\
&&\hspace*{80.38pt}{}
- \frac{16 \sigma_m(x) }{(2 \pi)^{1/2} \delta_m(x)}\exp\biggl\{-\frac
{(\delta_m/\sigma_m)^2}{8} \biggr\}\biggr),\nonumber
\end{eqnarray}
where the last inequality follows from Lemma
\ref{lm:boundRS_EPP_new}
[if $\delta_m(x) \rightarrow0$ and $m p_m \rightarrow1$
then for any $(y,x)$ there exists $M$ such that
$\forall m \geq M$, $C_{\delta_m(x)}(y) \cap A_0^m(x) = \varnothing$ and
the lemma applies].
Convergence of the bound in (\ref{eq:M4_lb}) to $f(y|x)$ a.s. $F$
is implied by a.s. positivity and continuity in $y$ of $f(y|x)$ and
conditions in (\ref{eq:cond_delta_sigma_h_x}).
%The uniform convergence in $x$ is not needed for point-wise
%convergence; however, it is used in deriving an integrable upper bound
%for DCT applicability.
The rest of the argument establishing point-wise convergence is the
same as for $\mathcal{M}_0$ [details are below (\ref
{eq:cond_delta_sigma_h})].

Next, let us derive an integrable upper bound for the DCT,
%
%e8.16 ###
%
\begin{eqnarray}\label{eq:RiemannSumEPPInBound}\hspace*{28pt}
p(y|x ,\mathcal{M}_4) &=&
\sum_{j=1}^m F(A_j^m(x)|x) \phi(y, \mu_j^m(x),\sigma_m(x))\nonumber\\
&&{}+F(A_0^m(x)|x) \phi(y, 0,\sigma_0(x)) \nonumber\\
&\geq&
[1-1_{A_0^m(x)}(y)] \nonumber\\
&&\hspace*{0pt}{}\times\inf_{z \in C_1(r(x),y,x)}
f(z|x)\nonumber\\[-8pt]\\[-8pt]
&&\hspace*{66pt}{}\times
\sum_{j\dvtx A_j^m(x) \subset C_1(r(x),y,x)} \lambda(A_j^m(x))\nonumber\\
&&\hspace*{161pt}{}\times \phi(y,
\mu
_j^m(x),\sigma_m(x))
\nonumber\\
&&{} + 1_{A_0^m(x)}(y) \cdot
\inf_{z \in C_0(r(x),y,x)} f(z|x) \cdot
\lambda(C_0(r(x),y,x))\nonumber\\
&&\hspace*{117.2pt}{}\times \phi(y, 0,\sigma_0(x)). \nonumber
\end{eqnarray}

Lemma \ref{lm:boundRS_EPP_new} and condition
(\ref{eq:cond_r_sigma_h_x})
imply that
the sum in (\ref{eq:RiemannSumEPPInBound}) is bounded
below by $1/2 - 1/4=1/4$ for all sufficiently large $m$.
Equation (\ref{eq:cond_sigma0_x}) implies
%
%e8.17 ###
%
\begin{eqnarray}\label{eq:EPPinqlty1}\quad
&&\log\max\biggl\{1,\frac{f(y|x)}{p(y|x ,\mathcal{M}_4)}\biggr\}
\nonumber\\
&&\qquad\leq
\log\max\biggl\{1,\frac{f(y|x) \cdot(r(x)/2)^{-1}}{\inf_{z \in C(r(x),y,x)}
f(z|x) \cdot\phi(y, 0,\sigma_0(x)) }\biggr\} \nonumber\\
&&\qquad\leq
\log\frac{1}{\phi(y, 0,\sigma_0(x)) (r(x)/2)}
\max\biggl\{\phi(y, 0,\sigma_0(x)) \bigl(r(x)/2\bigr),\\
&&\hspace*{184pt} \frac{f(y|x)}{\inf_{z \in
C(r(x),y,x)} f(z|x)}\biggr\} \nonumber\\
&&\qquad\leq
-\log[\phi(y, 0,\sigma_0(x)) r(x)/2 ] +
\log\frac{f(y|x)}{\inf_{z \in C(r(x),y,x)} f(z|x)}.\nonumber
\end{eqnarray}
Inequality (\ref{eq:EPPinqlty1}) follows by (\ref{eq:cond_sigma0_x}).
The first expression in (\ref{eq:EPPinqlty1})
is integrable by
Assumption \ref{assn:flexiblemu_case}, part \hyperlink
{assnitem:flexiblemu_case_2}{6}.
The second expression in (\ref{eq:EPPinqlty1}) is integrable by
Assumption \ref{assn:flexiblemu_case}, part \hyperlink
{assnitem:flexiblemu_case_3}{3}.
This completes the proof of the proposition.
\end{pf*}
\begin{pf*}{Proof of Corollary \protect\ref{crl:flex_mean_bddsupport}}
It suffices to show that Assumption \ref{assn:flexiblemu_case} is satisfied.
First, let us obtain a suitable $h_m$.
Note that
%
%e8.18 ###
%
\begin{equation}\qquad
\label{eq:p_mbd1}
p_m \geq\int_{A_j^m(x) \cap[a_1(x) , b_1(x)]} f(y|x)\,dy \geq\lambda
\bigl(A_j^m(x) \cap[a_1(x) , b_1(x)]\bigr) \underline{f}.
\end{equation}
Also,
%
%e8.19 ###
%
\begin{eqnarray}\label{eq:p_mbd2}
p_m &\geq& \int_{A_j^m(x) \cap[a(x),a_1(x)]} f(y|x)\,dy\nonumber\\
&\geq&\int_{A_j^m(x) \cap[a(x),a_1(x)]} \underline{f} \cdot[
y-a(x)]^n \,dy
\\
&\geq& (n+1)^{-1} \lambda\bigl(A_j^m(x) \cap[a(x),a_1(x)]\bigr)^{n+1}
\underline{f}\nonumber
\end{eqnarray}
%
%(the integral is bounded below by the area under $\underline{f} \cdot
%[ y-a(x)]^n$)
and similarly
$p_m \geq(n+1)^{-1} \lambda(A_j^m(x) \cap[b_1(x),b(x)])^{n+1}
\underline{f}$.
Combining this inequality with
(\ref{eq:p_mbd1}) and (\ref{eq:p_mbd2})
we get for all $x$ and $j$,
\begin{eqnarray*}
\lambda(A_j^m(x)) &\leq& \frac{p_m}{\underline{f}} + \frac{2 \cdot
(n+1)^{1/(n+1)} \cdot p_m^{1/(n+1)}}{\underline{f}^{1/(n+1)}}\\
&\leq&
\frac{7 p_m^{1/(n+1)}}{\underline{f}} = h_m.
\end{eqnarray*}
For $\sigma_m(x)=p_m^{1/4(n+1)}$ and $\delta_m(x)=p_m^{1/8(n+1)}$
conditions (\ref{eq:EPPconditions}), (\ref{eq:cond_delta_sigma_h_x})
and (\ref{eq:cond_r_sigma_h_x}) hold.

Next, let $C(r,y,x)=[y, y+r]$ if $y \in(a(x), a_1(x)+r/2)$,
$C(r,y,x)=[y-r/2, y+r/2]$ if $y \in[a_1(x)+r/2, b_1(x)-r/2]$
and
$C(r,y,x)=[y-r/2, y]$ if $y \in(b_1(x)-r/2, b(x))$.
By condition \hyperlink{crlitem:monotone}{4} of the corollary
$\inf_{z \in C(r(x),y,x)} f(z|x) = f(y|x)$ for $y \notin[a_1(x)+r/2,
b_1(x)-r/2]$.
For $y \in[a_1(x)+r/2, b_1(x)-r/2]$,
$\inf_{z \in C(r(x),y,x)} f(z|x) \geq\underline{f}$ and
\[
\int\log\frac{f(y|x)}{\inf_{z \in C(r(x),y,x)} f(z|x) } F(dy,dx)
\leq\log(\overline{f}/\underline{f}) < \infty.
\]
Condition \hyperlink{assnitem:flexiblemu_case_1}{2} and (\ref
{eq:cond_sigma0_x}) in Assumption \ref{assn:flexiblemu_case} are assumed
in the corollary.
Since
$a(x)$ and $b(x)$ are assumed to be square integrable,
the second moment of $y$ is finite, and
condition \hyperlink{assnitem:flexiblemu_case_2}{6} of Assumption \ref
{assn:flexiblemu_case} holds.
\end{pf*}
%
%Using an argument identical to the one in the proof of Corollary
%error in the following form
%& c_1 \cdot(\frac{1}{m} )^{\alpha_{\delta}}
%+c_2 \cdot(\frac{1}{m} )^{\alpha_{h} +2 \alpha_{\delta} - 3
%+c_3 \cdot(\frac{1}{m} )^{2 \alpha_{h} + \alpha_{\delta} -
%3 \alpha_{\sigma}} \\
%&+c_4 \cdot(\frac{1}{m} )^{\alpha_{\sigma} -\alpha_{\delta}}
%+c_5 \cdot(\frac{1}{m} )^{(q-2)(1-\alpha_{h})}, \nonumber
%where constants $c_1$, $c_2$, $c_3$, and $c_4$ do not depend on $m$.
%The largest common lower bound for powers of $1/m$ in
%%\[ \frac{1}{5+1/(q-2)}.\]
%
\begin{lemma} \label{lm:boundRSbyInt_ESP}
Define a hypercube $C_{\delta}(y)=\{\mu\in R^d\dvtx y_i \leq\mu_i \leq
y_i + \delta, i=1,\ldots,d\}$.
Let $A_1,\ldots,A_m$ be adjacent hypercubes with centers $\mu_j$ and
side length $h$ such that $C_{\delta}(y) \subset\bigcup_{j=1}^m A_j$ and
$\delta> 3 d^{1/2} h $. Define $J=\{j\dvtx A_j \subset C_{\delta}(y)\}
$. Then
\[
\sum_{j \in J}{\lambda(A_{j})\phi(y;\mu_{j},\sigma)}
\geq \int
_{C_{\delta}(y)}{\phi(\mu;y;\sigma)\,d\mu}-\frac{3 d^{3/2} \delta
^{d-1}h}{(2\pi)^{d/2}\sigma^{d}}.
\]
By symmetry, this result holds for any hypercube with vertex at $y$ and
side length~$\delta$. This implies that for hypercube $D_{\delta
}(y)=\{
x\dvtx y_i - \delta/2 \leq x_i \leq y_i + \delta/2, i=1,\ldots,d \}$,
\[
\sum_{j\dvtx A_j \subset
D_{\delta}(y)}{\lambda(A_{j})\phi(y;\mu_{j},\sigma)}
\geq
\int_{D_{\delta}(y)}{\phi(\mu;y;\sigma)\,d\mu}-2^d \frac{3 d^{3/2}
(\delta/2)^{d-1}h}{(2\pi)^{d/2}\sigma^{d}}
\]
as long as $D_{\delta}(y) \subset\bigcup_{j=1}^m A_j$ and $\delta> 6
d^{1/2} h$.
\end{lemma}
\begin{pf}
For $j \in J$ let $B_j = \{ x\dvtx\mu_{ji} \leq x_i \leq\mu_{ji}+h,
i=1,\ldots,d\}$ be a shifted and rotated version of $A_j$. Note that
$\mu_j = \arg\max_{\mu\in B_j} \phi(\mu;y;\sigma)$, and therefore
\begin{eqnarray*}
&&\sum_{j \in J}{\lambda(A_{j})\phi(y;\mu_{j},\sigma)}\\
&&\qquad=
\sum_{j \in J}{\lambda(B_{j})\phi(y;\mu_{j},\sigma)}
\geq
\int_{\bigcup_{j \in J} B_j}{\phi(\mu;y;\sigma)\,d\mu}\\
&&\qquad\geq
\int_{C_{\delta}(y)}{\phi(\mu;y;\sigma)\,d\mu}
-\int_{C_{\delta}(y) \setminus\bigcup_{j \in J} B_j}{\phi(\mu
;y;\sigma
)\,d\mu}.
\end{eqnarray*}
Since\vspace*{1pt}
$
\{x\dvtx\min_J \mu_{ji} \leq x_i \leq\max_J \mu_{ji}, i=1,\ldots,d \}
\subset C_{\delta}(y) \cap[\bigcup_J B_j]
$ and
$\max_{j \in J} \mu_{ji} - \min_{j \in J} \mu_{ji} \geq\delta- 3
d^{1/2} h$,
we get\vspace*{2pt}
$\lambda(C_{\delta}(y) \cap[\bigcup_J B_j]) \geq(\delta- 3 d^{1/2}
h)^d$ and
\begin{eqnarray*}
\lambda\biggl(C_{\delta}(y) \Bigm\backslash\biggl[\bigcup_J B_j\biggr]\biggr)
&=& \lambda(C_{\delta}(y)) - \lambda\biggl(C_{\delta}(y) \cap\biggl[\bigcup_{j \in
J} B_j\biggr]\biggr)
\nonumber\\
&\leq&
\delta^d - (\delta- 3 d^{1/2} h)^d \leq3 d^{3/2} h \delta^{d-1},
%3 h d \delta^{d-1}
\end{eqnarray*}
where the last inequality follows by induction.
Thus,
\begin{eqnarray*}
\int_{C_{\delta}(y) \setminus\bigcup_J B_j}{\phi(\mu;y;\sigma)\,d\mu}
&\leq&
\lambda\biggl(C_{\delta}(y) \Bigm\backslash\biggl[\bigcup_J B_j\biggr]\biggr) \frac{1}{(2 \pi
)^{d/2} \sigma^d} \\
&\leq&
\frac{3 d^{3/2} h \delta^{d-1}}{(2 \pi)^{d/2}\sigma^d}.
\end{eqnarray*}
\upqed\end{pf}
\begin{lemma} \label{lm:boundIntBy1}
Let $C_{\delta}(y)$ be a $d$-dimensional hypercube with center $y$ and
side length $\delta>0$. Then
\[
\int_{C_{\delta}(y)}{\phi(\mu;y;\sigma)\,d\mu} > 1 -\frac{8d\sigma
/\delta
}{(2\pi)^{1/2}} \exp\biggl\{-\frac{(\delta/\sigma)^2}{8} \biggr\}.
\]
Note that this inequality immediately implies that
for any sub-hypercube of $C_{\delta}(y)$, $\tilde{C}$, with vertex at
$y$ and side length $\delta/ 2$, for example,
$\tilde{C} = C_{\delta}(y) \cap[\mu\geq y]$,
\begin{eqnarray*}
\int_{\tilde{C}}{\phi(\mu;y;\sigma)\,d\mu} &=&
\frac{1}{2^{d}} \int_{C_{\delta}(y)}{\phi(\mu;y;\sigma)\,d\mu} \\
&>&
\frac{1}{2^{d}} -\frac{8d\sigma/\delta}{ 2^{d} (2\pi)^{1/2}}
\exp\biggl\{-\frac{(\delta/\sigma)^2}{8} \biggr\}.
\end{eqnarray*}
\end{lemma}
\begin{pf}
\begin{eqnarray*}
\int_{C_{\delta(y)}}{\phi(\mu;y;\sigma)\,d\mu} &=&
\int_{ \bigcap_{i=1}^d [|\mu_{i}| \leq\delta/2]} {\phi(\mu;0;\sigma
)\,d\mu}\\
&=& 1 - \int_{\bigcup_{i=1}^d [|\mu_{i}| \geq\delta/2]}
{\phi(\mu;0;\sigma)\,d\mu} \\
&\geq& 1 - \sum_{i=1}^{d}
\int_{|\mu_{i}| \geq\delta/2}
{\phi(\mu_{i};0;\sigma)\,d\mu_{i}}
\\
&=& 1- 2 d \int_{\delta/2 }^{\infty}{\phi(\mu_1;0;\sigma)\,d\mu_1}
\\
&>& 1 - \frac{2 d}{(2 \pi)^{1/2} \sigma} \int_{\delta/2 }^{\infty}{
\exp\biggl\{-\frac{0.5 (\delta/2) \mu_1 } {\sigma^2} \biggr\}\,d\mu_1}
\\
&=& 1 - \frac{2 d}{(2 \pi)^{1/2} \sigma}
\frac{- \sigma^2}{0.5 (\delta/2)}
\exp\{-0.5 (\delta/2) \mu_1 / \sigma^2 \}
|_{\delta/2}^{\infty}
\\
&=& 1 - \frac{8 d (\sigma/ \delta) }{(2 \pi)^{1/2} }
\exp\biggl\{-\frac{(\delta/\sigma)^2}{8} \biggr\}.
%> 1 - \frac{8 d \sigma}{(2 \pi)^{1/2} \delta}
\end{eqnarray*}
\upqed\end{pf}
\begin{lemma} \label{lm:boundRS_EPP_new}
Let $A_1,\ldots,A_m$ be a partition of an interval on $R$ such that
$\lambda(A_j) \leq h$ and $\mu_j \in A_j$. Assume $C_{\delta
}(y)=[y-\delta,y+\delta] \subset\cup A_j$ is an interval with center
$y$ and length $\delta$.
Then
\[
\sum_{j=1}^m \lambda\bigl(A_j \cap C_{\delta}(y)\bigr) \phi(y, \mu_j,\sigma
) \geq
1 -
\frac{ 6 h } {(2 \pi)^{1/2} \sigma}
% - \frac{8 \sigma}{(2 \pi)^{1/2} \delta}.
- \frac{8 (\sigma/ \delta) }{(2 \pi)^{1/2} }
\exp\biggl\{-\frac{(\delta/\sigma)^2}{8} \biggr\}.
\]
If $C_{\delta}(y)=[y-\delta,y]$ or $C_{\delta}(y)=[y,y+\delta]$ the
lower bound in the above expression should be divided by 2.
\end{lemma}
\begin{pf}
Let $J=\{j\dvtx A_j \cap C_{\delta}(y) \subset[y-\delta,y]\}$.
For any $j \in J$ and $\mu\in A_j \cap C_{\delta}(y)$, $\mu- h \leq
\mu_j$ as $\lambda(A_j) < h$ and $\mu_j \in A_j$, which implies
$\phi(y, \mu_j,\sigma) \geq\phi(y, \mu-h,\sigma)$. Therefore,
%
%e8.20 ###
%
\begin{equation}
\label{eq:sum_bd_shift}\qquad
\sum_{j \in J} \lambda\bigl(A_j \cap C_{\delta}(y)\bigr) \phi(y, \mu
_j,\sigma)
\geq
\int_{\bigcup_{j \in J} [A_j \cap C_{\delta}(y)]} \phi(y, \mu
-h,\sigma)\, d
\mu.
\end{equation}
Note next that
\begin{eqnarray*}
&& \int_{\bigcup_{j \in J} [A_j \cap C_{\delta}(y)]} \phi(y, \mu
-h,\sigma)
\,d \mu\\
&&\qquad\geq
\int_{y-\delta}^{y-h} \phi(y, \mu-h,\sigma)\, d \mu= \int
_{y-\delta
-h}^{y-2h} \phi(y, \mu,\sigma) \,d \mu
\\
&&\qquad =
\int_{y-\delta}^{y} \phi(y, \mu,\sigma)\, d \mu\\
&&\qquad\quad{} -
\int_{y-\delta-h}^{y-\delta} \phi(y, \mu,\sigma)\, d \mu-
\int_{y-2h}^{y} \phi(y, \mu,\sigma)\, d \mu
\\
&&\qquad \geq
\int_{y-\delta}^{y} \phi(y, \mu,\sigma) \,d \mu- \frac{3h}{(2 \pi
)^{1/2} \sigma}.
\end{eqnarray*}
By symmetry the same results can be obtained for $J=\{j\dvtx A_j \cap
C_{\delta}(y) \subset[y,y+\delta]\}$. Thus
\[
\sum_{j=1}^m \lambda\bigl(A_j \cap C_{\delta}(y)\bigr) \phi(y, \mu_j,\sigma)
\geq
\int_{y-\delta}^{y+\delta} \phi(y, \mu,\sigma) \,d \mu- 2 \frac{3h}{(2
\pi)^{1/2} \sigma}.
\]
The claim of the lemma follows by Lemma \ref{lm:boundIntBy1}.
\end{pf}
\end{appendix}

\section*{Acknowledgments}
The author is grateful to John Geweke and participants of seminars at
Princeton, SBIES 09 and SITE 09 for helpful discussions. I thank
Justinas Pelenis for pointing out shortcomings in several proofs. I
thank an associate editor and anonymous referees for useful
suggestions. All remaining errors are mine.

\printaddresses

\end{document}